%% file: main.tex
\documentclass[11 pt, draftcls, onecolumn]{IEEEtran}
\include{header}
\usepackage{mathtools}
\mathtoolsset{showonlyrefs=true}
\UseRawInputEncoding
\def \rh  {\hat r}
\begin{document}
%=======================================================================================================================
%\title{Popularity-Guaranteed Rebate Program}
\title{\vspace{-36pt}
{Sales-Based Rebate Design}\thanks{We are grateful to Dirk Bergemann, Matthew Jackson, Ben Golub, Ali Kakhbod, David Simchi-Levi, Alireza Tahbaz-Salehi,  and Rakesh Vohra for helpful comments. This work was supported by ARO MURI  W911NF-12-1-0509. \vspace{0.1cm}}
}

%Price Discrimination over Quality via Aggregate Rebate Programs

\author{Amir Ajorlou$^\dagger$
\and Ali Jadbabaie$^\dagger$\\
%\and David Simchi-Levi$\dagger$\\
%{\small{Working Draft}}
\thanks{$\dagger$Institute for Data, Systems, and Society, Massachusetts Institute of Technology (MIT), Cambridge, MA 02139, USA. E-mail: %\url{{ajorlou,jadbabai}@mit.edu.}}
% \url{{ajorlou,jadbabai}@mit.edu.}
}
%\thanks{This work was supported by ARO MURI  W911NF-12-1-0509.}
}
\date{September 2021}
%=======================================================================================================================
\maketitle
\thispagestyle{empty}
\normalsize
%%%%%%%%%%%%%%%%%%%%%%%%%%%%%%%%%%%%%%%%%%%%%%%%%%%%%%%%%%%%%%%%%%%%%%%%%%%%%%%%
\vspace{-46pt}
\begin{abstract}
\begingroup
%     \fontsize{8.9pt}{10pt}\selectfont

%Firms offering network products and services can enjoy a free gain in their profit resulted from the positive effect of additional users on the value of the product to others.
%Motivated by this observation,
%we pose the question as to whether it is beneficial to a firm selling a product with no (or weak) network effect to artificially inject externality?
We propose a novel family of sales-based rebate mechanisms that induce network effects in sales of products that do not exhibit such externalities. The proposed rebate mechanisms enable the seller of a product with uncertain quality to adjust the magnitude and sign of externalities in consumers' payoffs by conditioning the amount of the rebate on the sales volume. Using the machinery of global games and variational optimization techniques, we  analyze the revenue implications of such induced externalities in the form of rebate.  We identify optimal profitable designs while unraveling the main drivers of profit and further elucidate the main practical barriers associated with their implementation and show how these difficulties can be handled. The key insight of our rebate design is monetizing the strategic uncertainty in consumers' beliefs on others' valuations. The common externality induced in consumers' utilities as a sales-based rebate essentially enables the seller to elicit different prices at different valuations, given the heterogeneity of the beliefs on sales volume and hence on the rebate. Our analysis indicates that a mechanism that creates positive externalities will in fact reduce the profit because it lowers the expected prices at high valuations. On the other hand,  sellers can use a sales-based rebate mechanism that is decreasing with sales volume to incentivize purchase at lower valuations by providing higher expected rebates. Our work contributes to the literature  on technology-enabled features of digitized markets and demonstrates that  real-time sales/subscription data can lead to new revenue management methods.
\endgroup
\end{abstract}
\vspace{-12pt}
\begin{keywords}
\vspace{-9pt}
Network effects, strategic uncertainty, global games, group buying, variational optimization, rebates.
\end{keywords}

%%%%%%%%%%%%%%%%%%%%%%%%%%%%%%%%%%%%%%%%%%%%%%%%%%%%%%%%%%%%%%%%%%%%%%%%%%%%%%%
\section{Introduction}
\label{sec::intro}
Over the past two decades, the study of network effects and externalities has received a lot of attention both in microeconomics and in operations management. It is now well-understood and a part of conventional wisdom that goods and products that exhibit network effects or positive externalities can realize higher profits to sellers. With such strategic complementarities, consumer's utility from consumption of such goods and services is often higher than the intrinsic value of the product~\cite{CandoganBimpikisAsu.2012, kakhJad}.
%Early buyers improve the social value of the product and firm can charge a higher price when the product reaches certain level of popularity.
New technologies and innovations such as smartphone applications (e.g.,
WhatsApp, Signal, Telegram), online games (e.g., Warcraft), social
networking websites (e.g., Facebook, Twitter, Instagram,TickTok), and online dating
services (e.g., Zoosk, Match.com, OkCupid) are among many
examples of products with positive network effects.
The economic theory of network externalities has a long history, going back to   \cite{Farrell_RAND_85, Katz_JPE_86}. These works are followed by a series of papers on network games with strategic complements\footnote{Games of strategic complementarities are those in which the best response of each player is increasing in actions of others (\cite{Vives_90}).} (\cite{BallesterCalvo-ArmengolZenou.2006, Sundararajan_BEJTE_2007, GGJVY.2010}). More recently,
\cite{CandoganBimpikisAsu.2012, Bloch_GEB_2013, Harsha_working_2017, Wang_working_2013,kakhJad} study the effect of network externalities on optimal pricing and revenue management over networks in monopoly settings.

Motivated by the operational value of such network effects, it is natural to ask whether such effects can be induced in cases where they do not naturally exist.
% a natural question is whether it is beneficial for a firm selling a product with no (or weak) network effect to somehow create payoff externalities.
To answer this question, we introduce  a family of \textit{sales-based rebate mechanisms} as a means of inducing payoff externalities among consumers by conditioning the amount of the rebate on sales volume. 

Our setting consists of a seller and a continuum of buyers. Buyers purchase the product at a fixed exogenous market price and receive a rebate (according to the publicly announced structure by the seller) at the end of the sales period. 
% Buyers pay a fixed market price when making a purchase and are paid a reward at the end of the sales period. The reward is designed to be a function of the ex-post aggregate size of the buyers 
% Conditioning the amount of the rebate on the 
% (\emph{sales volume}) seller can induce externalities in utilities, thus resembling the network effect.
The proposed sales-based rebate mechanism enables us to formally model the effect of creating payoff externalities for products with no inherent network effect. 
% This enables us to mathematically model the effect of creating payoff externalities 
% for products with no inherent network effect. 
% This enables us to formulate the problem of optimally creating payoff externalities 
% for products with no inherent network effect.
We use this model to derive several key characteristics for the optimal design, and to investigate revenue implications of sales-based reward programs in general.
A key finding is that, while firms gain higher profits for products with inherent positive externalities, inducing such an effect by other means (e.g., a reward program) may indeed have an adverse impact on the profit.

Consumers in our model are differentiated by their intrinsic valuations of the product, parameterizing both the objective average quality of the product and subjective taste-related preferences of consumers. 
Both firm and consumers are uncertain about the average quality of the product, sharing a common prior belief on the quality.   

In this setting, each consumer privately observes her own valuation, but cannot separate the common quality component from her idiosyncratic taste component. Therefore, she has incomplete information
about the valuation of others and hence their purchase decisions.
The combination of heterogeneous tastes and uncertainty in the common quality induces a \textit{global game} with correlated private valuations among consumers\footnote{See, e.g., \cite{Pavan_RESTUD_19, Rossella_RAND_08} for other applications of global games in revenue management and pricing.} (\cite{MorrisShin.1998, MorrisShin.2003, CarlssonVanDamme.1993}). We assume that the firm does not possess any private information about the quality. 
This allows us to abstract away from the signaling implications of prices and reward programs, and focus instead on the social operations management aspects that stem from the heterogeneous yet strategic responses of consumers to uncertainty in purchase decision of others.
% and thus on the expected sales volume. 

The heterogeneity of the beliefs on purchase decisions of others (and hence on the sales volume) 
enables the seller to effectively price-discriminate, inducing different expected net prices at different valuations.
We can use this observation to show that, contrary to what one might believe at first,
introducing positive externality using a rebate program that is increasing with the sales volume may indeed reduce the seller's profit, since it will effectively result in lower expected prices at higher valuations.
%On the other hand,  for any decreasing non-constant reward function, the price can be adjusted in such a way that the joint price-reward program is profitable. Especially, we show that the optimal program can pay back a rebate as high as a full refund if the product fails to reach certain level
%of popularity.
%The main idea here is that a carefully designed reward program can induce higher net prices (in expectation) at higher valuations, unlike the fixed price case which charges all the buyers the same rate.

A carefully designed reward program can induce higher prices in expectation at higher valuations. 
Finding the optimal reward program, however, involves solving an infinite dimensional non-concave maximization problem with a continuum of constraints. While fully solving for the optimal solution in closed form is not  possible, we identify several key characteristics of the optimal reward program using variational optimization techniques. In particular, we show that the optimal reward program is a ``full-refund or nothing'' policy which
pays back the full price to buyers if the realized average quality lies in one of the finitely-many refund-eligible intervals.
The number of intervals, though finite, grows unboundedly as consumers' tastes become less diverse and valuations concentrate further around the common quality. 
Moreover, limiting the number of refund-eligible intervals to reduce design complexity can substantially degrade the performance of the optimal design, approaching that of the no-reward case in this regime.

These limitations along with challenges in implementing a non-monotone sales-based reward program, urge the need for implementable yet efficient alternatives to the optimal design. 
We propose and characterize two such alternatives by analytically solving for the optimal solution within two subspaces of sales-based reward functions:
one with a constraint on the reward spread  and the other with a constraint on its rate of change. Both reward programs turn out to be monotone, have simple structures, and perform provably-well when compared to the optimal solution;  the former coincides with the optimal solution when most of the uncertainty in consumers' valuations comes from the diversity in tastes, while the latter yields an asymptotically optimal profit when uncertainty in valuations is primarily rooted in the common quality.

\subsection{Literature Review}
\label{subsec::litreview}

This work is closely related to several areas in the revenue management and pricing literature, including group-buying and quantity discounts. Our proposed reward program, however, offers a different source of revenue compared to  group-buying schemes.% which for a while became quite popular.
~In a group-buying scheme, buyers can receive a discount if they simultaneously  purchase the product  as a  group. In 2010, \emph{Groupon}, a major player in group-buying industry, was named  the fastest-growing company in the history of the Web by Forbes (\cite{Forbes_2010}). Despite their stunning early rise, the industry has experienced a downfall over the past few years: \emph{LivingSocial} (Groupon's main competitor), once valued at \$6 billion, was acquired by Groupon for \$0 in 2016 (\cite{Wharton_2017}). Groupon's stock value dropped from a high of more than \$28 in 2011 to below the \$1 mark in 2020, triggering the clock on the first step of getting delisted from \textit{Nasdaq}. The company had to execute a 1-for-20 reverse stock split in Summer 2020. \emph{Amazon Local}, one of Groupon's competitors offering similar packages, closed down in 2015 (\cite{GeekWire_2015}).
Given the ups and downs of the group-buying industry, the pros and cons of their business model has been debated, discussed, and dissected, yet the future of group-buying platforms is an uncertain one.%yet to be seen.

%The benefit of group-buying schemes is mainly attributed to its reliance on the ``economies of
Benefits of group-buying strategies are often pointed out as ``economies of
networking'' and ``economies of scale'' in the business press (\cite{Forbes_2012}). Along the same line,
\cite{Xie_MS_2011} suggest that the key advantage of group-buying lies in fostering Word-of-Mouth: it incentivizes the experts to act as ``sales agents'' and to promote the product to novice customers through interpersonal influences. Such strategies appear more suitable for relatively unknown firms, as was shown in \cite{Edelman_Marketing_2016}.
In a related work, \cite{Pinar_Omega_2016} study group-buying mechanisms by explicitly accounting for both the utility from shopping together with one's social circle\footnote{See \cite{Bristol_Retailing_2004} and the references therein for the influence of the peers on spending more on shopping.} as well as the inconvenience cost due to the wait time, and show that the former usually outweighs the latter.
\cite{Kauffman_JMIS_2001} find evidence of the positive
externality effect on customer bids using customer data from %\href{http://www.mobshop.com/}
{MobShop.com}.
Selling in large groups is also advantageous in situations involving scale economies (e.g., in restaurant's industry) as large quantities reduce the marginal cost (\cite{ Monahan_MS_1984, Kohli_MS_1989}).

In this work, however, our focus is on products with fixed market size and marginal cost.
This enables us to single out the operational value of reward programs coming from their direct effect on the utilities of the firm and consumers in absence of scale economies, while setting aside second order effects such as market expansion via Word-of-Mouth and interpersonal influences. In this setting, our results suggest that discounting the price as a function of the size of buyers is no more a profitable strategy.

Another approach to group-buying and threshold discounting is to view them as means of dealing with demand uncertainty. \cite{Aron_MS_2003} and \cite{Zhang_MS_2015}
use threshold discounting to find the operative demand regime in a scenario where the seller is uncertain about the demand. Unlike our setting, however,
demand parameters are assumed to be fully known to buyers.
Treating the entire market as a single player with unknown type, where the type determines the operative demand regime, their results can be closely related to the seminal work of \cite{Maskin_RAND_1984} which studies optimal quantity discounting of a seller in face of a buyer with uncertain type.\footnote{There are still substantial differences in the assumptions on the type distribution in these works. As a result, while quantity discounts in \cite{Maskin_RAND_1984} are everywhere optimal, threshold discounting strategies  may not always outperform posted fixed prices as noted by the authors in \cite{Aron_MS_2003} and \cite{Zhang_MS_2015}.} Somewhat closer to our work is that of \cite{Karan_INSEAD_2016}, where demand uncertainty is present at both ends. A seller with capacity constraints
uses threshold discounting to both signal the market size to buyers and to condition offering the product during the ``slow'' season on the market size, hence reducing the supply-demand mismatch.\footnote{\cite{Cachon_MS_2004} also uses quantity discounts to encourage early season purchases to reduce the
risk due to demand uncertainty.} As noted by the authors, however, this strategy can potentially reduce the profit if the seller has no capacity constraint, as assumed in our work.

Another related body of work is the literature on \emph{referral reward programs}, where the seller uses monetary rewards to motivate existing
buyers to spread product information thus expanding the market (\cite{Libai_Marketing_2001, Aral_MS_2011, Ilan_MS_2016, Jackson_working_2017}).\footnote{Other word of mouth marketing strategies include creating \emph{buzz} using promotions  and frequent \emph{zero-pricing} (see \cite{CampbellMayzlinShin.2013}, and \cite{Amir_MS_2017} and references therein).}
Although very similar in nature, group-buying has
the advantage of stimulating a larger scale of social
interaction as it requires information sharing before any transaction takes place (see \cite{Xie_MS_2011} for a detailed comparison of group-buying and referral reward programs). As stated before, by considering a fixed market size we set aside the second order beneficial marketing effects of the reward programs, including the market expansion via social interactions, in our analysis and fully elaborate on the operational value of such programs resulted from their direct effect on the utilities of the seller and buyers.

In summary, 
while the operational value of network effects is mostly attributed to  
its effectiveness in fostering word of mouth, social influence, and scale economies, the value in creating externalities here comes from monetizing the strategic uncertainty (\cite{MorrisShin_2002}) - that is, uncertainty concerning the purchase decision of customers, as well as  beliefs (and beliefs about the beliefs) on purchase decisions of other customers.
This value can be realized by inducing a properly-designed common externality component into consumers' payoffs as a ``sales-based rebate'', which enables the seller to induce different prices at different valuations given the heterogeneity of the beliefs on sales volume and hence on the rebate.
While firms achieve a growth in their profit for products with an inherent positive network effect, we articulate that inducing such an effect using a sales-based reward program may reduce the profit, as it will effectively induce lower expected prices at higher valuations.
A seller, on the other hand, may benefit from a properly designed decreasing rebate function (though may not be optimal), exploiting its capacity to induce lower expected prices at lower valuations.
Our work thus complements the overwhelming literature on the benefits of positive network effects by putting spotlight on the operational value of creating negative payoff externalities. 
Nevertheless, incentive programs such as group-buying and referral rewards
can be still beneficial due to their effectiveness in fostering word of mouth and social influence, scale economies, and reducing supply-demand mismatch under capacity constraint
in situations discussed in the literature of group-buying and referral reward programs.

%We use a simple stylized model to highlight the core ideas of our proposed sales-based rebate program. Recent advances in technology have made it possible to keep track of the users' statistics specially when it comes to digital goods and services, a feature that was not available a few years back. Our work aims at supporting a theory for a new generation of rebate programs that takes advantage of these new features.\footnote{The closest implementation that we have found is that of ``Guaranteed prize pool'' poker tournaments in online poker (e.g., \href{http://www.pokerstarsnj.com/poker/tournaments/sunday-tournaments/}{PokerStars} and \href{http://www.fulltilt.com/poker/tournaments/types}{FullTilt}), where the platform (also known as the \emph{house}) guarantees a certain number of participants for the tourney. Each player pays a fee to register for the tourney which goes to the prize pool. If the prize pool falls short of the promised size, then the rest is on the house.}
%

Along with their analytical complexity and operational challenges, technology-driven markets bring a series of useful features that were not previously available.
%This suggests the need to revisit/improve some of the traditional revenue management strategies as they may no longer be efficient in the larger space of implementable strategies to date.
Our work aims at developing a theory for a new generation of rebate programs that takes advantage of these new features.
We use a simple, stylized model to highlight the core ideas of our proposed sales-based rebate mechanism,
which requires keeping track of the number of users of a product -a feature easy to implement today (at least for digital goods and services). Such a feature would be quite out of reach a few years ago.\footnote{The closest implementation that we have found is that of ``Guaranteed prize pool'' poker tournaments in online poker (e.g., PokerStars and FullTilt, %\href{http://www.pokerstarsnj.com/poker/tournaments/sunday-tournaments/}{PokerStars} and \href{http://www.fulltilt.com/poker/tournaments/types}{FullTilt}),
where the platform (also known as the \emph{house}) guarantees a certain number of participants for the tourney. Each player pays a fee to register for the tourney which goes to the prize pool. If the prize pool falls short of the promised size, then the rest is on the house.}

%seller cannot control informational environment (as in bayesian persuasion) but can control their payoffs.

%``However, the implementation of such contracts requires the
%principal to cross-check the reports by (potentially) all players to detect an individual’s deviation.
%This can be daunting when the network size is large.''

%Useful phrases:
%Our objective is to use a simple model to illustrate the ...
%foster, stimulate, boost, enhance,
%differences in tastes
%positive externality is to enjoy but not to create unless for second order effects such as fostering word of mouth or scale economics.
%In contrast, creating negative externality can effectively discriminate the price over quality.
%our work departs from this literature

%%%%%%%%%%%%%%%%%%%%%%%%%%%%%%%%%%%%%%%%%%%%%%%%%%%%%%%%%%%%%%%%%%%%%%%%%%%%%%%%
\section{Model}
\label{sec::model}
We consider a firm selling an indivisible product to a unit-mass continuum of consumers indexed by $i\in[0,1]$ at an exogenous market price $p$.\footnote{We elaborate on this point  later in the section.}  
Consumers are differentiated by their intrinsic valuation of the product. 
Product valuation of consumer $i$ is given by $v_i=v+\epsilon_i$, where the common component $v$ captures the objective/average quality of the product, and $\epsilon_i$ 
represents the subjective/taste-related preferences consumers may have for the product. 

Both the firm and consumers are uncertain about the average quality, sharing a common Gaussian prior belief $v\sim N(\theta,\sigma_\theta^2)$ on the quality of the product.   
Each consumer privately observes her own valuation, but cannot separate the common quality component from her idiosyncratic taste component, which we assume has a normal distribution $\epsilon_i\sim N(0,\sigma_\ep^2)$.

The two sources of uncertainty in a consumer's valuation, that is, i) the uncertainty in the common quality (quantified with $\sigma_\theta$), and ii) an {idiosyncratic uncertainty} resulting from the diversity of tastes (parametrized with $\sigma_\ep$),  yields a total uncertainty of variance $\sigma^2=\sigma_\ep^2+\sigma_\theta^2$, and hence an ex-ante distribution $v_i\sim N(\theta,\sigma^2)$ for consumers' valuations.

%We may interchangeably refer to the quality of a product as its \textit{quality}.
%Firm and consumers share a common prior $v\sim N(\theta,\sigma_\theta^2)$ on the quality of the product.

When there is no reward program, 
an agent with valuation $v_i$ makes a purchase ($a_i=1$) if the utility of purchase, given by $u_i=v_i-p$, is higher than the utility of not buying ($a_i=0$), which is normalized to zero. This results in a \textit {sales volume} of size
$\bar a(v)=\Prob [v_i>p|v]=\Phi(\frac{v-p}{\sigma_\ep})$ for a given realization $v$ of the average quality, where $\Phi(\cdot)$ denotes the CDF of the standard normal distribution.
This yields an ex-ante expected sales volume of size
\begin{align}
\E_v[\bar a(v)]=\E_v[\Phi(\frac{v-p}{\sigma_\ep})]=\Phi(\frac{\theta-p}{\sigma}).
\end{align}
%where $\sigma^2=\sigma_\ep^2+\sigma_\theta^2$ is the variance of the total uncertainty in consumers' valuations.
The expected profit of the seller is thus,
\begin{equation}
\bar \Pi(p)=p\E_v[\bar a(v)]=p\Phi(\frac{\theta-p}{\sigma}).\footnote{The marginal cost is normalized to zero.}
\end{equation}

% \begin{Assumption}
% \label{ass::marketprice}
% The market price $p$ is set so as to maximize the ex-ante expected profit of the seller, and is given by the unique solution of
% \begin{align}
% \frac{p}{\sigma}=\frac{\Phi(\frac{\theta-p}{\sigma})}{\phi(\frac{\theta-p}{\sigma})},
% \end{align}
% where $\sigma^2=\sigma_\ep^2+\sigma_\theta^2$ is the variance of the ex-ante total uncertainty in consumers' valuations.
% % We assume an optimally set fixed price for the product, in order to purely isolate the effect of the reward program.}
% \end{Assumption}

To model the effect of sales-based reward programs on the profit, 
% Our aim here is to study the effect of reward programs on the profit of the seller.
% The timing of the actions of agents in the model is as follows:
we assume that the firm announces a reward program  $r:[0,1]\rightarrow[0,p]$ when launching the product at an exogenous market price $p$.
% at the market price $p$, which is assumed to be optimally set as in Assumption~\ref{ass::marketprice}.
Consumers pay the price $p$ to the firm when they buy the product knowing that they will receive a reward valued at $r(\bar a)$ at the end of the sales period, where $\bar a$ is the ex-post sales volume.
%We also assume that the aggregate reward $R(\bar a)=r(\bar a)\bar a$ is increasing with the sales volume, implying that firm does not make withdrawals from the reward pool.

%\debate{$\theta$ can be regarded as the average valuation in the market.}
Upon observing the price and reward program,
% and given the private valuations and the common public prior on the quality,
agents simultaneously decide whether to purchase or not. 
% The payoff  for non-buyers ($a_i=0$) is normalized to zero, as before. 
The offered reward program induces a new component into the payoff of a purchase ($a_i=1$):
\begin{equation}
\label{eq::payoff}
u_i=v_i+r(\bar a)-p.
\end{equation}
The offered reward program conditions the utility that a consumer derives from a purchase on purchase decisions of others, thereby resembling the network effect. 
Consumers take actions maximizing their expected payoffs, speculating on the expected sales volume. The information available to customer $i$ at the time of making her purchase decision are: i) her own valuation $v_i$, ii) the common prior on the average/common quality, and iii) the price and the announced reward function. We can express the utility-maximizing purchase decisions as
\begin{equation}
\label{eq::actions}
a_i={\bf 1}\{\E_v[u_i|v_i,p,r(\cdot)]>0\}.
\end{equation}
The utility of the firm offering the reward program $r(\cdot)$ is then
\begin{equation}
\label{eq::firmutility}
\Pi(p,r(\cdot))=(p-r(\bar a))\bar a,
\end{equation}
where $\bar a$ is the sales volume resulted from the  purchase strategies of the consumers.
% and the realized quality $v$.
%The expected utility of the firm is thus $\E_v[\Pi(p,r(\cdot),v)|\theta]=\E_v[(p-r(\bar a(v)))\bar a(v)|\theta]$.

Our aim here is to use the above setting to gain insight on how a seller can use a (carefully designed) sales-based reward program in order to strategically control the induced network effect of its product. 
Specifically, we are interested in characterizing the optimal design (maximizing the expected profit) and profitable reward programs in general;
we call a reward program \textit{profitable} if it yields an expected profit higher than the profit of the no-reward case. 
% A natural question which arises is whether there exists a profitable reward program at all. 
% One can think of a reward program as a bet between the firm and each consumer whose outcome depends on the realization of the quality. Therefore, a profitable reward program may seem out of reach, given the informational ``edge" of the consumers over the firm.\footnote{Besides the common prior on the quality $v$, each consumer's valuation provides her with a noisy private observation of $v$.}
\subsection{The case with no uncertainty in average quality ($\sigma_\theta=0$)}
% \section{Profitability of Aggregate Reward Programs}
A key feature of our model is its ability to account for 
the heterogeneity of consumers' beliefs on purchase decisions of others, and hence on the sales volume.
Absent such heterogeneity, consumers can fully coordinate on the expected sales volume and consequently the amount of the reward they will receive. This effectively degrades the reward program to a mere price discount. 

To illustrate this further, consider an orthodox setting where 
the average quality $v$ is perfectly known (i.e., $\sigma_\theta=0$ and $v=\theta$). 
Consumers will therefore share a common belief $v_i\sim N(\theta, \sigma_\theta^2)$
on distribution of valuations. This leads to a common prediction on the expected sales volume $\bar a$, and hence the reward $r(\bar a)$, at any equilibrium induced by the offered reward program $r(\cdot)$. Consumers will subsequently internalize this, adjusting the price and making a purchase if and only if $v_i>c=p-r(\bar a)$. This makes the offered reward program effectively equivalent to a discount lowering the price from $p$ to $c=p-r(\bar a)$.

In order to focus on the nontrivial benefits of a sales-based reward program, we make the following optimality assumption on the market price $p$. 

\begin{Assumption}
\label{ass::marketprice}
The market price $p$ is set so as to maximize the ex-ante expected profit of the seller, and is given by the unique solution of
\begin{align}
\frac{p}{\sigma}=\frac{\Phi(\frac{\theta-p}{\sigma})}{\phi(\frac{\theta-p}{\sigma})},
\end{align}
where $\sigma^2=\sigma_\ep^2+\sigma_\theta^2$ is the variance of the ex-ante total uncertainty in consumers' valuations.\footnote{ $\phi(\cdot)$ denotes the PDF of the standard normal distribution.}
\end{Assumption}

%This means that no reward program can beat the profit resulted from optimally setting the price with no reward.
%
%\begin{Proposition}
%\label{prop::knowntheta}
%Consider the case where the public signal $v$ fully reveals the quality of the product, that is $\sigma_\theta=0$. Then, given price $p$, an increasing reward program $r(\cdot)$ satisfying Assumption~\ref{ass::uniqueness}, and quality of the product $\theta$, we have
%\begin{equation}
%\Pi(p,r(\cdot),\theta)\leq\Pi(p_c^*,0,\theta),
%\end{equation}
%where $p_c^*=\argmax_{p}p\Phi(\frac{\theta-p}{\sigma_\epsilon})$.
%\end{Proposition}

The combination of heterogeneous tastes and uncertainty in the average quality induces a {global game} with correlated private valuations among consumers. 
The first step in analyzing the sales-based reward programs is to characterize the equilibria of the subgame among the consumers induced by the offered reward program.
% strategic responses of consumers to uncertainty in purchase decision of others.

\subsection{Monotone Bayes-Nash equilibria of the consumers' subgame}
The negligible effect of individual consumers on the aggregate action in continuum models makes the Bayes Nash equilibria of the game symmetric.
We specifically turn our attention to equilibria within the class of \textit{monotone} or \textit{threshold} strategies.
A symmetric, monotone strategy with threshold $c$ is of the form $a_i={\bf 1}\{v_i>c\}$. For such a strategy, a consumer makes a purchase if and only if her private valuation is above the threshold $c$.
% The existence of such equilibria for the consumers' subgame can be guaranteed by imposing certain conditions on the reward program, as we will see later in Lemma~\ref{lemma::rewardspread} and ~\ref{lemma::rewardROC}.
We next elaborate on how to characterize the threshold equilibria induced by a given reward program, by analyzing the utility-maximizing decision making of consumers following such a purchase strategy. 
%We first derive conditions ensuring the existence of such equilibria for the consumers' subgame.

% For a monotone purchase strategy $a_i={\bf 1}\{v_i>c\}$,
% the aggregate size of the buyers for a realization $v$ of the quality is $\bar a=\Phi(\frac{v-c}{\sigma_\ep})$.
Observing her private valuation $v_i=v+\epsilon_i$, consumer $i$ updates her belief on the common quality from the prior $v\sim N(\theta,\sigma_v^2)$ to $v|v_i\sim N(\tau v_i+(1-\tau)\theta,\sigma_v^2)$, where $\tau=\frac{\sigma_\theta^2}{\sigma_\ep^2+\sigma_\theta^2}$ and $\sigma_v^2=\frac{\sigma_\ep^2\sigma_\theta^2}{\sigma_\ep^2+\sigma_\theta^2}$. 
A monotone purchase strategy $a_i={\bf 1}\{v_i>c\}$,
yields a sales volume of size $\bar a=\Phi(\frac{v-c}{\sigma_\ep})$ for a realization $v$ of the quality. To see this, note that
\begin{align}
    \bar a ~=~& {\rm Prob}\big(v_i> c| i\text{ is a random customer}\big) = {\rm Prob}(v+\ep_i> c)\nonumber\\
    =~&{\rm Prob}(\ep_i> c-v) = {\rm Prob}(\ep_i< v-c)=\Phi(\frac{v-c}{\sigma_\ep}).
\end{align}
Putting this together with her updated belief on the quality and the offered reward program, she can then speculate on the expected payoff yield from a purchase:
\begin{equation}
\label{eq::expectedpayoff}
\E_{v|v_i}[u_i]= v_i-p+\E_{v|v_i}[r(\Phi(\frac{v-c}{\sigma_\ep}))].
\end{equation}
Given an offered reward program $r(\cdot)$, a monotone purchase strategy $a_i={\bf 1}\{v_i>c\}$ is thus an equilibrium strategy if and only if
\begin{gather}
\label{eq:threshEQ}
\E_{v|v_i}[r(\Phi(\frac{v-c}{\sigma_\ep}))]+v_i-p\geq0 \text{ for } v_i\geq c,\\
\E_{v|v_i}[r(\Phi(\frac{v-c}{\sigma_\ep}))]+v_i-p\leq0 \text{ for } v_i\leq c.
\end{gather}
The above constraints are to ensure that consumers' purchase decisions are 
utility-maximizing responses to the purchase strategy $a_i={\bf 1}\{v_i>c\}$, thus making it an equilibrium strategy.
The consumer with valuation $c$, or so-called the \textit{cutoff}, 
is  indifferent to making a purchase or not. We can write the \textit{indifference equation}  $\E_v[u_i|v_i=c]=0$ as $c-p+r_c=0$,
% \begin{align}
% \label{eq::cutoff}
% c-p+r_c=0,
% \end{align}
where $r_c=\E_{v|v_i=c}[r(\Phi(\frac{v-c}{\sigma_\ep}))]$ is the reward expected at the cutoff.

Each buyer is charged a
{``net price''} of $p-r(\Phi(\frac{v-c}{\sigma_\ep}))$.
The expected utility of the firm thus becomes:
\begin{align}
\label{eq::firmexpectedutility}
\E_v[\Pi(p,r(\cdot))]&=\E_v[(p-r(\Phi(\frac{v-c}{\sigma_\ep}))) \Phi(\frac{v-c}{\sigma_\ep})]\nonumber\\
&=\E_v[(p-\hat r(v)) \Phi(\frac{v-c}{\sigma_\ep})].
%&=c\E_v[\Phi(\frac{v-c}{\sigma_\ep})]+
%\E_v[(r_c-r(\Phi(\frac{v-c}{\sigma_\ep})))\Phi(\frac{v-c}{\sigma_\ep})]\nonumber\\
%&=\E_v[\Pi(c,0)]+
%\E_{v|c}[r(\Phi(\frac{v-c}{\sigma_\ep}))]\E_v[\Phi(\frac{v-c}{\sigma_\ep})]- \E_v[r(\Phi(\frac{v-c}{\sigma_\ep}))\Phi(\frac{v-c}{\sigma_\ep})].
\end{align}

The continuum of constraints in \eqref{eq:threshEQ} identifies the feasible space of reward programs for which the induced subgame among the consumers admits a threshold equilibrium purchase strategy, while \eqref{eq::firmexpectedutility} shows the profit the firm can expect from appending the reward program to its price. This formulation enables us to gain insight into the profit implications of sales-based reward programs by viewing it in a variational optimization framework (\cite{Ito_2008, Optimization_Luenberger, Clarke_2013}).  

It is more convenient to work with reward as a function of the quality $v$. Given the one-to-one map between the quality and sales volume ($\bar a=\Phi(\frac{v-c}{\sigma_\ep})$) for threshold strategies, we define and henceforth work with $\hat r(v)=r(\Phi(\frac{v-c}{\sigma_\ep}))$ in the rest of the paper.

%Another useful representation is to write the expected utility of the firm as the sum (integral) of the expected profit made at each valuation $v_i$ weighted by its density among all the realizations of $v$ and $\epsilon_i$ resulting in $v_i$:
%that is writing the profit by conditioning on the realization of the valuation $v$ versus, conditioning on the valuation $v_i$ over all realizations of $v$, i.e., systemic versus idiosyncratic view of the profit.
%

\section{Optimal Sales-Based Reward Program}
The optimal sales-based reward program is the solution to the following infinite-dimensional optimization problem:
% We can formulate the problem of finding the optimal reward program as
\begin{gather}
\label{eq::optimal}
\underset{\substack{\hat r\in L^\infty(\mathbb{R};[0,p]),c\in\mathbb{R}}}{\text{maximize }} \E_v[(p-\rh(v)) \Phi(\frac{v-c}{\sigma_\ep})],\\%\E_v[\Pi(p,\rh(v))],\\
\text{subject to:}\\
\E_{v|v_i}[\rh(v)]+v_i-p\geq0 \text{ for } v_i\geq c,\\
\E_{v|v_i}[\rh(v)]+v_i-p\leq0 \text{ for } v_i\leq c,
\end{gather}
%where $L_1(\mathbb{R};[0,p],d\normcdf{v-\theta}{\sigma_\theta})$ is the space of $\Phi$-measurable functions from $\mathbb{R}$ to $[0,p]$.
where $L^\infty(\mathbb{R};[0,p])$ is the space of bounded measurable functions taking values in $[0,p]$. It is to be noted that the above formulation allows for price adjustment (if necessary) while designing the reward program, as long as the offered price is below the market price $p$ and the reward is capped with the offered price.   
More precisely, any joint price-reward program pair $(\hat r, \hat p)$ with $0\leq \hat p\leq p$ and $\hat r \in L^\infty(\mathbb{R};[0,\hat p])$ can be effectively implemented, as it is equivalent to offering an adjusted reward program $\hat r+p-\hat p$ at the market price $p$.

In this section, we use the above formulation to identify several key characteristics of the optimal sales-based reward program. We highlight the key steps of our approach in dealing with 
the optimization problem in \eqref{eq::optimal} and refer the readers to the appendix for the details.

\subsection{Optimality of ``full-refund or nothing'' reward programs}
\label{subsec::optimality}
A common approach in dealing with constrained optimization problems as in \eqref{eq::optimal} is based on Lagrange multiplier theory (see, e.g., \cite{Ito_2008, Optimization_Luenberger, Clarke_2013}). 
Existence of Lagrange multipliers for this problem can be established using a regularity condition that basically requires the linearized approximation of the constraint space around the optimal solution to have a feasible interior point:\footnote{see, e.g., Definition 1.5 in \cite{Ito_2008} for the explicit statement of the regularity used here.}
Denote the optimal cutoff associated with \eqref{eq::optimal} by $c^*$. One can easily verify that for the fixed reward $\rh(v)= p-c^*$ and $c=c^*$ all the inequalities are strict and the indifference equation is satisfied. Since the only nonlinearity in \eqref{eq::optimal} (indifference equation for the cutoff) is kept unchanged, this verifies the aforementioned regularity condition.
Let the Lagrangian be
\begin{align}
\label{eq::Lagrangian}
L(\rh,c,\lambda)=\E_v[(p-\rh(v)) \Phi(\frac{v-c}{\sigma_\ep})]+\int_{\mathbb{R}}(c+x+\E_{v|v_i=c+x}[\hat r(v)]-p)d\lambda(x),
\end{align}
where $\lambda\in BV(\mathbb{R})$ is an upper-semicontinuous function with bounded variation, decreasing for $x<0$ and increasing for $x\geq0$.
Denote the optimal solution of \eqref{eq::optimal} with $(\rh^*,c^*)$.\footnote{Let $X=L^1(\mathbb{R},d\normcdf{v-\theta}{\sigma_\theta})$. Existence of a global maximizer follows from the
weak$^*$ compactness of the closed unit sphere in $X^*=L^\infty(\mathbb{R})$ and the continuity of the expected profit in $(\rh,c)\in L^\infty(\mathbb{R};[0,p])\times [0,p]$ (see the proof of Theorem~\ref{theorem::optimal-characteristics} for details).}
The complementary slackness property requires
\begin{align}
\label{eq::CS}
\int_{\mathbb{R}}(c^*+x+\E_{v|v_i=c^*+x}[\hat r^*(v)]-p)d\lambda(x)=0.
\end{align}
Noting the nonnegativity of the integrand, this implies that $\lambda(x)$ can only change value when $c^*+x+\E_{v|v_i=c^*+x}[\hat r^*(v)]-p=0$, that is, when a consumer with valuation $v_i=c+x$ is indifferent between making a purchase or not. It is easy to see that the indifference equation can only admit finitely many solutions $x$: Any such solution should lie in $[0,p]$; on the other hand, the payoff of consumers is an analytic function of $x$ (since $\E_{v|v_i}[r(v)]$ is analytic), and hence can only admit finite number of zeros in  $[0,p]$. Denote the set of such values of $v_i=c+x$ with $x\in\mathcal X$. Then, noting that $d\lambda(x)=0$ for $x\notin \mathcal X$, we can write the Lagrangian as
\begin{align}
\label{eq::Lagrangian2}
L(\rh,c,\lambda)=\E_v[(p-\rh(v)) \Phi(\frac{v-c}{\sigma_\ep})]+\sum_{x\in\mathcal X}\hat\lambda(x)(c+x+\E_{v|v_i=c+x}[\hat r(v)]-p),
\end{align}
where $\hat\lambda(x)=\lambda(x)-\lambda(x^-)$, with $\hat\lambda(x)<0$ for $x<0$ and $\hat\lambda(x)>0$ for $x>0$.
Expanding \eqref{eq::Lagrangian2}, we can rewrite the Lagrangian as
\begin{align}
\label{eq::Lagrangian3}
L(\rh,c,\lambda)=&p\normcdf{\theta-c}{\sigma}+\sum_{x\in\mathcal X}\hat\lambda(x)(c+x-p)\nonumber\\
&-\int_{\mathbb{R}}\rh(v)\left(\normcdf{v-c}{\sigma_\ep}\normpdf{\theta-v}{\sigma_\theta}-\sum_{x\in\mathcal X}\hat\lambda(x)\normpdf{\mu_c+\tau x-v}{\sigma_v}\right)dv.
\end{align}
Optimality of $(\hat r^*,c^*)$ requires:
\begin{gather}
\label{eq::FOC1}
\frac{\partial}{\partial c}L(\rh^*,c,\lambda)_{|c=c^*}=0,\\
L(\rh^*,c^*,\lambda)\geq L(\rh,c^*,\lambda),\text{ for all }\rh\in L^\infty(\mathbb{R};[0,p]).
\end{gather}
Optimal reward program is hence a ``full-refund or nothing'' policy of the form
\begin{equation}
\label{eq::optimalr}
\rh^*(v)=p\times{\bf 1}\{g(v)<0\},
\end{equation}
where
\begin{align}
\label{eq::g(v)}
g(v)=\normcdf{v-c^*}{\sigma_\ep}\normpdf{\theta-v}{\sigma_\theta}-\sum_{x\in\mathcal X}\hat\lambda(x)\normpdf{\mu_c+\tau x-v}{\sigma_v}.
\end{align}
This function encapsulates the cost-benefit analysis for changes in the value of the reward at quality $v$: the burden on the profit of the seller is determined by i) how likely is for $v$ to be the realized quality ($\normpdf{\theta-v}{\sigma_\theta}$), and ii) the sales volume at this quality ($\normcdf{v-c^*}{\sigma_\ep}$). 
Lagrange multiplier $\hat \lambda(x)$ measures the marginal cost of violating the constraint at valuation $v_i=c^*+x$ with a marginal change in the utility of the consumer. This makes it nonzero only at indifferent valuations, positive above the threshold $c^*$ and negative below it. 
It then needs to be adjusted by the sensitivity of the expected utility of the respective consumer to the value of reward at quality $v$, and is captured by $\normpdf{\mu_c+\tau x-v}{\sigma_v}$.

Linearity of both the seller's profit and the constraints in the reward makes the optimal solution a ``full-refund or nothing'' policy; if having reward at quality $v$ is found to be costly overall (i.e., $g(v)>0$) then $\hat r^*(v)=0$, and if profitable (an overall negative cost, i.e., $g(v)<0$) then $\hat r^*(v)=p$. It is to be noted that $g$ can have multiple zeros in general, resulting in a non-monotone  optimal reward program.

%and $\eta(v|v_i)=\frac{\phi(\frac{v-\mu_i}{\sigma_v})}{\sigma_v}$ denotes the (posterior) belief of a consumer with valuation $v_i$ on the quality.
Using the first order condition for the optimal threshold $c^*$ in \eqref{eq::FOC1} we can find
\begin{align}
\hat\lambda(0)=\frac{c^*\phi(\theta-c^*)}{\frac{\partial}{\partial c}h(0,c)_{|c=c^*}}>0,
\end{align}
where $h(x,c)=c+x+\E_{v|v_i=c+x}[\hat r(v)]-p$.\footnote{Note that for a primal feasible reward program, $\frac{\partial}{\partial x}h(x,c^*)=0$ for every nonzero $x\in\mathcal{X}$. This, in turn, yields $\frac{\partial}{\partial c}h(x,c)_{|c=c^*}=0$ for every nonzero $x\in\mathcal{X}$.}
%\begin{align}
%\label{eq::gasymptotic1}
%\lim_{v\to\infty} \frac{g(v)}{\frac{\phi(\frac{\theta-v}{\sigma_\theta})}{\sigma_\theta}}=1.
%\end{align}
%Similarly,
%\begin{align}
%\label{eq::gasymptotic2}
%\lim_{v\to-\infty} \frac{g(v)}{\normpdf{\ubar\mu_x-v}{\sigma_v}}=-\hat\lambda(\ubar x),
%\end{align}
%where $\ubar x=\min\{x|x\in\mathcal X\}$.
Asymptotic behavior of the function $g(v)$ can be seen from
\begin{align}
\label{eq::gasymptotic}
\lim_{v\to\infty} \frac{g(v)}{\frac{\phi(\frac{\theta-v}{\sigma_\theta})}{\sigma_\theta}}=1,\qquad\lim_{v\to-\infty} \frac{g(v)}{\normpdf{\mu_{\ubar x}-v}{\sigma_v}}=-\hat\lambda(\ubar x),
\end{align}
where $\ubar x=\min\{x|x\in\mathcal X\}$.
The optimal reward program is of the form
\begin{align}
\rh^*(v)=p\times\sum_{j=1}^l {\bf 1}\{w_L^j\leq v\leq w_H^j\},
\end{align}
with $w_L^1=-\infty$ if and only if $\mathbb{R^-}\cap \mathcal X\neq\emptyset$.
This means that the optimal reward program is a {\it full-refund or nothing} policy which refunds the full price to a buyer if the realized quality $v$ falls in one of the $l$ refund-eligible intervals $[w_L^j,w_H^j]$, $j=1,\ldots,l$.

Given the non-concave infinite dimensional nature of the maximization problem in \eqref{eq::optimal}, it is quite tempting to wish for the optimal solution to have only a few number of full-refund intervals $l$ or at least, for such strategies to achieve a significant fraction of the optimal profit.  
Verifying whether or not such a desirable property holds requires establishing a connection between the performance of a full-refund or nothing policy and the number of its full-refund intervals.

\subsection{Performance of {\it full-refund or nothing} policies with $l$ refund intervals}
\label{subsec::performanceF}
Let $\mathcal V_{\rm refund}=\cup_{j=1}^l [w_L^j,w_H^j]$ denote the set of qualities eligible for a full refund for a full-refund or nothing policy with $l$ refund intervals.
Consider consumer $i$ with valuation $v_i\in [c,p]$ and let $\mu_i=\tau v_i+(1-\tau)\theta\in[\mu_c,\mu_p]$ be  the quality expected by this consumer.
Notice that, these are the subset of buyers who have made the purchase relying on some nonzero (positive) expected reward, since their valuations alone are below the price while their valuations are above the purchasing cutoff. 
We can bound the reward expected by a consumer expecting quality $\mu_i$ based on their distance from $\mathcal V_{\rm refund}$.
Fix a radius $\delta>0$ and let ${\rm B}(\mathcal V_{\rm refund}, \delta)$ be the set of all qualities within $\delta$-vicinity of $\mathcal V_{\rm refund}$, that is,
\begin{align*}
{\rm B}(\mathcal V_{\rm refund}, \delta):=\{\mu\in \mathbb{R}| d(\mu, \mathcal V_{\rm refund})\leq \delta\}. 
\end{align*}
We claim that for a consumer $i$ with $\mu_i\notin {\rm B}(\mathcal V_{\rm refund}, \delta)$, we have $\E_{v|v_i}[\hat r(v)]<2p\normcdf{-\delta}{\sigma_v}$. 
To see this, i) recall that the consumer $i$'s belief on quality is  
$v|v_i\sim N(\mu_i, \sigma_v^2)$, a normal distribution centered at $\mu_i$, and ii) notice that no reward is paid at qualities within $\delta$-vicinity of the center of this distribution. This results in an upper bound on the reward equal to the tail probability in normal distribution at a distance $\delta$ away from the center times the full refund price, that is $2p\normcdf{-\delta}{\sigma_v}$.

This implies that, the quality expected by any buyer $i$ with valuation  $v_i\in[c,p-2p\normcdf{-\delta}{\sigma_v}]$ should lie within ${\rm B}(\mathcal V_{\rm refund}, \delta)$, or equivalently, $[\mu_c,\mu_p-2\tau p\normcdf{-\delta}{\sigma_v}]\subseteq {\rm B}(\mathcal V_{\rm refund}, \delta)$. As a result,
\begin{align}
\label{eq::refundeligible1}
|[\mu_c,\mu_p]\setminus{\rm B}(\mathcal V_{\rm refund}, \delta)|<2p\tau\normcdf{-\delta}{\sigma_v},
\end{align}
thus limiting the mass of the qualities $\mu_i\in[\mu_c,\mu_p]$ that could lie outside the (one-dimensional) ball of radius $\delta$ around ${\rm B}(\mathcal V_{\rm refund}, \delta)$.
Observe that when $\delta=0$, \eqref{eq::refundeligible1}  becomes trivial since $\mu_p-\mu_c = \tau(p-c)<\tau p$, not providing any information on refund qualities. A useful property here is that the mass of qualities between $\mathcal V_{\rm refund}$ and ${\rm B}(\mathcal V_{\rm refund}, \delta)$ is upper-bounded by $2l\delta$, that is $|{\rm B}(\mathcal V_{\rm refund}, \delta)\setminus \mathcal V_{\rm refund}|\leq 2l\delta$. This simply follows from the fact that 
$\mathcal V_{\rm refund}$ is the union of $l$ intervals. Putting this together with \eqref{eq::refundeligible1}, we get
\begin{align}
\label{eq::refundeligible2}
|[\mu_c,\mu_p]\setminus\mathcal V_{\rm refund}|<2p\tau\normcdf{-\delta}{\sigma_v}+2l\delta,
\end{align}
for any choice of $\delta>0$. This provides us with a very useful upper-bound on the mass of qualities in $[\mu_c,\mu_p]$ that are not eligible for a full-refund.  
Notice that, these are the only values of the quality within $[\mu_c,\mu_p]$ that create profit for the seller, and how the number of refund-eligible intervals $l$ imposes a constraint on the mass of such qualities. 

We next use \eqref{eq::refundeligible2} to gain insight on the effect of the number of the refund-eligible intervals on the profit of the seller.  
The seller fully refunds the buyers if the realized quality $v\in \mathcal V_{\rm refund}$, and charges each buyer a net price equal to $p$ otherwise. Combining this with \eqref{eq::refundeligible2}, we can come up with the following upper bound on the profit of the seller:\footnote{See the proof of Theorem~\ref{theorem::optimal-characteristics} for details.}
\begin{align}
\label{eq::profuboundL}
\E_v[\Pi(p,\hat r(v))]<\frac{p}{\sqrt{2\pi}\sigma_\theta}(2p\tau\normcdf{-\delta}{\sigma_v}+2l\delta)
+ p\int_{\mathbb{R}\setminus[\mu_c,\mu_p]}\normcdf{v-c}{\sigma_\ep}\normpdf{\theta-v}{\sigma_\theta}dv
\end{align}
Fix the total uncertainty in consumers' valuations $\sigma^2=\sigma_\theta^2+\sigma_\ep^2$.
When $\sigma_\epsilon\to0$ the second term in the above upper bound will have an asymptotic value of $p\normcdf{\theta-p}{\sigma}$.
% , that is the profit of the no reward case. 
It then follows that
\begin{align}
\liminf_{\sigma_\epsilon\to0}\E_v[\Pi(p,\hat r(v))]\leq p\normcdf{\theta-p}{\sigma}+\liminf_{\sigma_\epsilon\to0}\frac{p}{\sqrt{2\pi}\sigma}(2p\normcdf{-\delta}{\sigma_\epsilon}+2l\delta).
\end{align}
Let us choose $\delta = \sqrt{\frac{\sigma_\epsilon}{l}}$. Then, 
\begin{align}
\liminf_{\sigma_\epsilon\to0}\E_v[\Pi(p,\hat r(v))]&\leq p\normcdf{\theta-p}{\sigma}+
\liminf_{\sigma_\epsilon\to0}{\frac{p}{\sqrt{2\pi}\sigma}(2p\normcdf{-1}{\sqrt{l\sigma_\epsilon}}+2\sqrt{l\sigma_\epsilon })}.
\end{align}
For fixed $l$, this implies that
\begin{align}
\liminf_{\sigma_\epsilon\to0}\E_v[\Pi(p,\hat r(v))]&\leq p\normcdf{\theta-p}{\sigma},
\end{align}
where  $p\normcdf{\theta-p}{\sigma}$ is the expected profit of the no reward case. Also, on noticing that for the optimal reward program, 
$\liminf_{\sigma_\epsilon\to0}\E_v[\Pi(p,\rh^*(v))]>\E_v[\Pi(p,0)]$, we find that, $\liminf_{\sigma_\epsilon\to0}{l\sigma_\epsilon}>0$. That is, $\frac{1}{\sigma_\epsilon}=\mathcal O(l)$, or equivalently $l=\Omega(\frac{1}{\sigma_\epsilon})$. Next theorem summarizes the main findings of this section. 

\begin{Theorem}
\label{theorem::optimal-characteristics}
The optimal reward program $\rh^*(\cdot)$ cast as the solution of the optimization problem in \eqref{eq::optimal}
is a ``full-refund or nothing'' policy of the form
\begin{align}
\rh^*(v)=p\times\sum_{j=1}^l {\bf 1}\{w_L^j\leq v\leq w_H^j\},
\end{align}
where the intervals $[w_L^j,w_H^j]$, $j=1,\ldots,l$ identify the qualities eligible for a refund. Moreover, fixing the total uncertainty in consumers' valuations $\sigma^2=\sigma_\theta^2+\sigma_\ep^2$, then  as $\sigma_\ep\to 0$,
\begin{itemize}%[leftmargin=*]
\item[i)]\noindent the number of refund-eligible intervals $l$ in the optimal reward program tends to infinity. More precisely, $\frac{1}{\sigma_\epsilon}=\mathcal O(l)$, or equivalently $l=\Omega(\frac{1}{\sigma_\epsilon})$. 
\item[ii)]\noindent for any fixed $l$, the expected profit resulted from the optimal reward program with $l$ full-refund intervals approaches that of the no reward case.
\end{itemize}
\end{Theorem}
\noindent\textit{Proof.} {See the appendix.}$\hfill\blacksquare$

\section{Monotone Sales-Based Reward Programs}
The possibility of having too many refund-eligible intervals in the optimal sales-based reward program and that limiting the number of such intervals may substantially degrade its performance, along with challenges in implementing it (because of its non-monotone structure) calls for easy-to-implement yet efficient alternatives. 

In this section, we propose two such alternatives by analytically solving for the optimal solution within two specific subspaces of the reward programs (Section~\ref{subsec::SC} and \ref{subsec::RC}), and compare their performances with that of the optimal design (Section~\ref{subsec::performanceA}). Both solutions turn out to be monotone functions of the sales volume. Monotonicity is clearly an appealing property when it comes to implementing a sales-based reward program. 
A monotone sales-based reward program is either an increasing function of the sales volume (and so of the quality\footnote{Recall the one-to-one map between the quality and sales volume ($\bar a=\Phi(\frac{v-c}{\sigma_\ep})$).}) thus resembling a positive network effect, or a decreasing function of the sales volume hence inducing negative externalities in consumers' payoffs.
A relevant question here is whether inducing externalities via monotone reward programs is indeed profitable?

We expressed the expected utility of the seller in \eqref{eq::firmexpectedutility} by conditioning on the realized quality of the product. We can come up with a useful alternative representation for the expected profit by (instead) conditioning on consumers' valuations.
Recall that $v_i\sim N(\theta,\sigma^2)$, where $\sigma^2=\sigma_\ep^2+\sigma_\theta^2$ is the total uncertainty in valuations, and that the net price expected by a consumer with valuation $v_i$ is $p-\E_{v|v_i}[\hat r(v)]$. Finally, let us use the indifference equation $c-p+r_c=0$, where $r_c = \E_{v|v_i=c}[\hat r(v)]$ is the reward expected at the cutoff, to substitute for price $p$ to get 
% $c+\E_{v|v_i=c}[\hat r(v)]-\E_{v|v_i}[\hat r(v)]$ as the net price expected by a consumer with valuation $v_i$. we obtain
\begin{align}
\label{eq::EPi}
\E_v[\Pi(p,\rh(v))]&=\int_{c}^{\infty}(p-\E_{v|v_i}[\hat r(v)])\psi(v_i)dv_i\nonumber\\
&=\int_{c}^{\infty}(c+\E_{v|v_i=c}[\hat r(v)]-\E_{v|v_i}[\hat r(v)])\psi(v_i)dv_i\nonumber\\
&=\E_v[\Pi(c,0)]+\int_{c}^{\infty}(\E_{v|v_i=c}[\hat r(v)]-\E_{v|v_i}[\hat r(v)])\psi(v_i)dv_i,
%\E_{v|c}[\hat r(v)]\E_v[\Phi(\frac{v-c}{\sigma_\ep})]- \E_v[\hat r(v)\Phi(\frac{v-c}{\sigma_\ep})].
\end{align}
where  $\psi(v_i)=\frac{\phi(\frac{\theta-v_i}{\sigma})}{\sigma}$ is the ex-ante pdf of consumers' valuations, and $\E_v[\Pi(c,0)]$ is the expected utility of the firm offering a price $c$ with no reward.
%Instead of optimizing over $(p,r(\cdot))$, we find the pair of effective price-reward program $(c,r(\cdot))$ maximizing the firm's profit given in \eqref{eq::firmexpectedutility} and then use the indifference equation  \eqref{eq::cutoff} to obtain the corresponding optimal price $p$.
The following result follows immediately.
\begin{Lemma}
\vspace{-4pt}
\label{lemma::necessarycond}
Suppose the consumers follow an equilibrium strategy of the form $a_i={\bf 1}\{v_i>c\}$, given the price-reward program pair $(p,\rh(\cdot))$. Then, $\E_v[\Pi(p,\rh(\cdot))]>\E_v[\Pi(c,0)]$ if and only if the ex-ante expected reward paid per purchase is less than the reward expected at the cutoff, that is,
\begin{equation}
\vspace{-10pt}
\label{eq::necessarycond}
\E_v[\hat r(v)|v_i=c]>\E_v[\hat r(v)|v_i\geq c],
\end{equation}
where
\begin{equation}
\vspace{-9pt}
\E_v[\hat r(v)|v_i\geq c]=\frac{\int_{c}^{\infty}\E_{v|v_i}[\hat r(v)]{\psi(v_i)}dv_i}{\int_{c}^{\infty}\psi(v_i)dv_i}.
\end{equation}
\end{Lemma}
\noindent\textit{Proof.} {See the appendix.}~$\hfill\blacksquare$

% The seller extracts the full surplus from the cutoff agents since the expected net price at cutoff is the same as their valuation $c$ (recall the indifference equation $c=p-r_c$). 
Writing the expected profit of the firm as in \eqref{eq::EPi}, we can observe that the extra surplus extracted from buyers with $v_i>c$ via the reward program is the difference in the reward paid at valuation $v_i$ and that paid at the cutoff (i.e., $\E_{v|v_i=c}[\rh(v)]-\E_{v|v_i}[\rh(v)]$). It is then clear that a reward program that pays higher reward at higher valuations cannot outperform an offered price $c$ with no reward, and hence cannot be profitable.

We can show the monotonicity of expected rewards for monotone reward programs (see the proof of Theorem~\ref{theorem::monotonereward}).
In particular, we can show that for any increasing non-constant reward program $\rh(\cdot)$ the expected rewards  are strictly increasing with consumers' valuations.
As a result, no increasing reward function can be profitable, since it clearly violates the condition in \eqref{eq::necessarycond}.
% , and hence cannot be profitable.
\begin{Theorem}
\label{theorem::monotonereward}
There exists no profitable increasing reward program.
\end{Theorem}

\noindent\textit{Proof.} {See the appendix.}~$\hfill\blacksquare$

This result shows that although an inherent positive network effect yields a significant boost in profit margins, inducing it via a sales-based reward program for a product with no (or weak) network effect is indeed harmful to the profit. For an increasing reward program, buyers with higher valuations are charged a lower net price  in expectation, an overtly non-profitable course of action.
A decreasing reward program, on the other hand, charges a higher expected price at valuations above the cutoff.
The offered reward program is hence profitable if the extra surplus extracted from the buyers via the reward program surpasses the
marginal loss in profit resulted from moving the cutoff away from the optimally set price $p$.

A primary source of complexity in dealing with the optimal sales-based reward program as formulated in \eqref{eq::optimal} is the continuum of constraints characterizing the feasible space of reward programs; that is, the reward programs for which the induced subgame among the consumers admits a threshold equilibrium purchase strategy.  
One simple idea to reduce this complexity is to restrict the solution to a subspace of the rewards with the property that  
(assuming any threshold purchase strategy for consumers)  the resulting expected purchase utilities will be monotone in consumer valuations.\footnote{This is referred to as the \textit{single-crossing property} in mechanism design literature.}
Restricting the reward program to such a subspace, the set of constraints in \eqref{eq:threshEQ} then boils down to a single constraint, that is the indifference equation for the cutoff.  
Below, we identify and analyze two such subspaces: one by imposing a constraint on the reward spread, and the other by imposing a constraint on its rate of change. 
% The optimal solutions in both cases turn out to be decreasing functions of the sales volume. 
% can ensure the monotonicity of the expected payoffs (in valuations) by either imposing conditions on the spread of the reward (as in Lemma~\ref{lemma::rewardspread}), or its rate of change (as in Lemma~\ref{lemma::rewardROC}).

%We next  present conditions under which $\frac{\partial}{\partial v_i}\E_v[u_i|v_i]\geq0$, that is, the expected payoff from a purchase is higher for a consumer with higher valuation.
\subsection{Spread-constrained sales-based reward programs}
\label{subsec::SC}

\begin{Lemma}
\label{lemma::rewardspread}
Let $r_{\min}$ and $r_{\max}$ denote the minimum and maximum reward paid to a buyer under the reward program $\rh(\cdot)$. Then, the expected payoff of adoption, assuming a monotone symmetric strategy for the consumers, is increasing with their valuations if

\begin{equation}
\label{ass::rewardspread1}
r_{\max}-r_{\min}\leq \sqrt{2\pi}\frac{\sigma_\ep\sigma}{\sigma_\theta}.\footnote{We can similarly come up with conditions ensuring the uniqueness of such an equilibrium. However, as we show in the proof of Theorem~\ref{theorem::optimal-BSpread}, the uniqueness condition will be trivial for the class of decreasing reward programs (that are of our special interest).}
\end{equation}
\end{Lemma}
\noindent\textit{Proof.} {See the appendix.}$\hfill\blacksquare$

%\begin{Assumption}
%\label{ass::rewardspread}
%$r_{\max}-r_{\min}\leq r_M$, where $r_M=\sqrt{2\pi}\sigma_\ep\sqrt{1+(\frac{\sigma_\ep}{\sigma_\theta})^2}$.
%\end{Assumption}

With this assumption on the reward range, the purchase strategy $a_i={\bf 1}\{v_i>c\}$ is an equilibrium strategy if and only if $c-p+\E_{v|v_i=c}[\rh(v)]=0$.
The optimal reward program under this assumption is  the solution to the following (infinite dimensional) optimization problem:
\begin{gather}
\label{eq::LPspread}
\underset{\hat r:\mathbb{R}\to[0,p], c\in\mathbb{R}}{\text{maximize }} \E_v[(p-\rh(v)) \Phi(\frac{v-c}{\sigma_\ep})],\\
\text{subject to:}\\
r_{\max}-r_{\min}\leq\frac{\sqrt{2\pi}\sigma_\ep\sigma}{\sigma_\theta},\\
c+\E_{v|v_i=c}[\rh(v)]=p.
\end{gather}
Exploiting the log-concavity of the normal distribution, we can analytically solve for the optimal reward program, as characterized in the next theorem.

\begin{Theorem}
\label{theorem::optimal-BSpread}
The optimal reward program $\rh^*_{SC}(\cdot)$ with the spread constraint
$r_{\max}-r_{\min}\leq \sqrt{2\pi}\frac{\sigma_\ep\sigma}{\sigma_\theta}$
is of the form $\rh^*_{SC}(v)=r_M\times{\bf 1}\{v\leq v_{c}\}$, where
$r_M=\min(p,\sqrt{2\pi}\frac{\sigma_\ep\sigma}{\sigma_\theta})$, and $\sigma=\sqrt{\sigma_\ep^2+\sigma_\theta^2}$,
 and $(c,v_c)$ is the solution of
\begin{align}
\label{eq::optimalreward}
p=&c+r_M\normcdf{v_c-\mu_c}{\sigma_v},\nonumber\\
\frac{\Phi(\frac{v_c-c}{\sigma_\ep})}{\frac{\phi(\frac{v_c-c}{\sigma_\ep})}{\sigma_\ep}}
=&\frac{c}{1-\tau r_M\normpdf{v_c-\mu_c}{\sigma_v}},
\end{align}
where $\sigma_v=\frac{\sigma_\ep\sigma_\theta}{\sqrt{\sigma_\ep^2+\sigma_\theta^2}}$, $\tau=\frac{\sigma_\theta^2}{\sigma_\ep^2+\sigma_\theta^2}$, and $\mu_c=\tau c+(1-\tau)\theta$.
\end{Theorem}

\noindent\textit{Proof.} {See the appendix.}~$\hfill\blacksquare$

According to this theorem, the optimal spread-constrained reward program
pays the buyers a fixed rebate with a value up to the full price if the realized quality (or the corresponding sales volume as for implementation purposes) falls below certain threshold $v_c$, and no rebate otherwise. Besides its simple structure, another useful property of the optimal spread-constrained reward program is the uniqueness of the equilibrium threshold strategy of the consumers' subgame under this reward program. This is a very useful property, since with multiple subgame equilibria seller would then need to use some other measure to speculate on whether the equilibrium strategy given by the optimal solution is the one capturing the purchase behavior of the consumers. 

The downside  is that the bound $\frac{\sqrt{2\pi}\sigma_\ep\sigma}{\sigma_\theta}$ on the reward spread in \eqref{ass::rewardspread1} is only affected by the diversity of tastes and the uncertainty in common quality, and is  independent from the expected quality of the product $\theta$ or, equivalently, the price $p$.
This may result in a reward which is quite insignificant compared to the paid price.
Consequently, the optimal reward program may perform poorly, especially in regimes where uncertainty in valuations is primarily rooted in the common quality and not the diversity in tastes ($\frac{\sigma_\ep}{\sigma}$ is small). 

\subsection{Rate-constrained sales-based reward programs}
\label{subsec::RC}
% The bound $\frac{\sqrt{2\pi}\sigma_\ep\sigma}{\sigma_\theta}$ on the reward spread in \eqref{ass::rewardspread1} is only affected by the diversity of tastes and uncertainty in average quality, and not the expected quality of the product $\theta$ or, equivalently, the price $p$.
% This may result in a reward which is quite insignificant compared to the paid price.
% Consequently, the optimal reward program may perform poorly, especially in regimes where
% the uncertainty in valuations is primarily due to the uncertainty in average quality and not the diversity in tastes ($\frac{\sigma_\ep}{\sigma}$ is small). 
The poor performance of the optimal spread-constrained reward program characterized in Theorem~\ref{theorem::optimal-BSpread} in regimes where
the uncertainty in valuations mainly comes from the uncertainty in the common quality
motivates searching for a reward program with a better performance in such regimes. To accommodate this, we look into another subspace of the reward programs which would still yield monotone expected utilities for buyers.

\begin{Lemma}
\label{lemma::rewardROC}
Let $\rh:\mathbb{R}\to[0,p]$ be a piece-wise continuously differentiable function with
\be
\label{eq::derivativebound}
|\frac{d}{dv}\rh(v)|\leq\frac{1}{\tau},
\ee
where $\tau=\frac{\sigma_\theta^2}{\sigma_\ep^2+\sigma_\theta^2}$.
Then, the expected payoff of adoption, assuming a monotone symmetric threshold strategy for the consumers, is increasing with their valuations.
\end{Lemma}
\noindent\textit{Proof.} {See the appendix.}$\hfill\blacksquare$

We can similarly characterize the optimal reward program within this subspace of reward programs using generalized Lagrange multipliers.

\begin{Theorem}
\label{theorem::optimallinear}
The optimal reward program $\rh^*_{RC}(\cdot)$ satisfying the rate of change constraint $|\frac{d}{dv}\rh(v)|\leq\frac{1}{\tau}$
is of the form
\begin{align}
\label{eq::roptimal}
\hat r^*_{RC}(v)=\begin{cases}
p, &v<v_L\\
p-\frac{1}{\tau}(v-v_L), &v_L\leq v<v_H\\
0, &v_H\leq v,
\end{cases}
\end{align}
where $(c,v_L,v_H)$ is the solution of
\begin{gather}
\label{eq::optimallinear}
\frac{\int_{v_L}^{v_H}{\normcdf{w-c}{\sigma_\ep}\normpdf{\theta-w}{\sigma_\theta}dw}}{\normcdf{v_H-\mu_c}{\sigma_v}-\normcdf{v_L-\mu_c}{\sigma_v}}=
\frac{c\phi(\theta-c)}{1+\normcdf{v_L-\mu_c}{\sigma_v}-\normcdf{v_H-\mu_c}{\sigma_v}},\\
v_H-v_L=\tau p,\\
p=c+\frac{\sigma_v}{\tau}\left((\frac{v_H-\mu_c}{\sigma_v})\normcdf{v_H-\mu_c}{\sigma_v}-(\frac{v_L-\mu_c}{\sigma_v})\normcdf{v_L-\mu_c}{\sigma_v}+
\normpdf{\frac{v_H-\mu_c}{\sigma_v}}{}-\normpdf{\frac{v_L-\mu_c}{\sigma_v}}{}\right),
\end{gather}
where $\sigma_v=\frac{\sigma_\ep\sigma_\theta}{\sqrt{\sigma_\ep^2+\sigma_\theta^2}}$, $\tau=\frac{\sigma_\theta^2}{\sigma_\ep^2+\sigma_\theta^2}$, and $\mu_c=\tau c+(1-\tau)\theta$.
\end{Theorem}
\noindent\textit{Proof.} {See the appendix.}$\hfill\blacksquare$

The above theorem states that, the optimal rate-constrained reward program is a decreasing function of the quality (or, equivalently sales volume) which fully refunds the customers if the realized quality falls below certain threshold $v_L$. The paid amount reduces at a constant rate, eventually reaching zero at some quality $v_H$ after which no rebate is paid.   
This rebate program specially proves efficient when uncertainty in valuations is primarily rooted in the common quality. 
In the extreme case when ($\sigma_\epsilon\to 0$), we can use the above theorem to show that $(c, v_L, v_H)\to (0, 0, p)$, implying that any customer with a positive valuation (though mostly concentrated around the realized quality) is incentivized to make a purchase. To see this, notice that as $\sigma_\epsilon\to 0$, 
customer $i$'s belief on the quality ($v|v_i\sim N(\tau v_i+(1-\tau)\theta,\sigma_v^2)$ concentrates around $v_i$, yielding an expected net price of $p-\hat r^*_{RC}(v=v_i)$ at valuation $v_i$. 
The rebate program with $(v_L, v_H)=(0, p)$ pays a rebate $\hat r^*_{RC}(v)=p-v$ when $v\in [0, p]$ and 0 otherwise, resulting in an expected net price equal to $v_i$ at valuation $v_i$ for $v_i\in [0, p]$, and a net price equal to $p$ for $v_i>p$.
This enables the seller to incentivize any customer with a positive valuation to purchase the product in the extreme case where $\sigma_\epsilon\to 0$.
% while fully extracting the surplus from those with valuation $v_i\in [0, p]$ and the full price $p$ from those with valuation above $p$. 

\subsection{Performance analysis}
\label{subsec::performanceA}
Unlike the optimal reward program formulated in \eqref{eq::optimal}, which is the solution of an infinite dimensional non-concave maximization problem with a continuum of constraints, the reward programs characterized in Theorem~\ref{theorem::optimal-BSpread} and \ref{theorem::optimallinear} have a simple structure, can be characterized analytically, and are easy to implement as they are both decreasing functions of the sales volume.
To evaluate the performance of these reward programs, however, we need  sufficiently tight upper bounds on the expected profit of the optimal sales-based reward program (that is, the solution to \eqref{eq::optimal}).

A simple yet useful observation is that the expected profit extracted from consumers with valuation $v_i$ is capped with both the price $p$ and $v_i$. This leads to the following upper bound on the expected profit:
\begin{align}
\label{eq::profubound1}
\Pi^H_1=&p\normcdf{\theta-p}{\sigma}+\int_{0}^{p}x\normpdf{\theta-x}{\sigma}dx\nonumber\\
=& p\normcdf{\theta-p}{\sigma}+ \theta(\normcdf{\theta}{\sigma}-\normcdf{\theta-p}{\sigma})+\sigma(\normpdf{\frac{\theta}{\sigma}}{}-\normpdf{\frac{\theta-p}{\sigma}}{}).
\end{align}

We can also obtain an upper bound on the optimal expected profit by solving a relaxed version of \eqref{eq::optimal}, in which we only keep the constraint corresponding to the indifference equation.\footnote{Alternatively, we could use weak duality.} %$\lambda(x)=\lambda_0\times{\bf 1}\{x\leq0\}$ in the Lagrangian in \eqref{eq::Lagrangian}.

\begin{Lemma}
\label{lemma::Profupperbound2}
Let $(c,v_c)$ be the solution of
\begin{align}
\label{eq::optimalreward}
c=&p\normcdf{\mu_c-v_c}{\sigma_v},\nonumber\\
\frac{\Phi(\frac{v_c-c}{\sigma_\ep})}{\frac{\phi(\frac{v_c-c}{\sigma_\ep})}{\sigma_\ep}}
=&\frac{c}{1-\tau p\normpdf{v_c-\mu_c}{\sigma_v}},
\end{align}
where $\sigma_v=\frac{\sigma_\ep\sigma_\theta}{\sqrt{\sigma_\ep^2+\sigma_\theta^2}}$, $\tau=\frac{\sigma_\theta^2}{\sigma_\ep^2+\sigma_\theta^2}$, and $\mu_c=\tau c+(1-\tau)\theta$.
Then the expected profit under the optimal reward program (that is, the solution to \eqref{eq::optimal}) cannot exceed the expected profit resulted from the
reward function $\rh(v)=p\times{\bf 1}\{v< v_{c}\}$ and cutoff $c$, given by
\begin{align}
\Pi^H_2=p\int_{v_c}^{\infty}\normcdf{v-c}{\sigma_\ep}\normpdf{\theta-v}{\sigma_\theta}dv.
\end{align}
\end{Lemma}
\noindent\textit{Proof.} {See the appendix.}$\hfill\blacksquare$

\begin{figure}
  \centering
  \includegraphics[scale=.55]{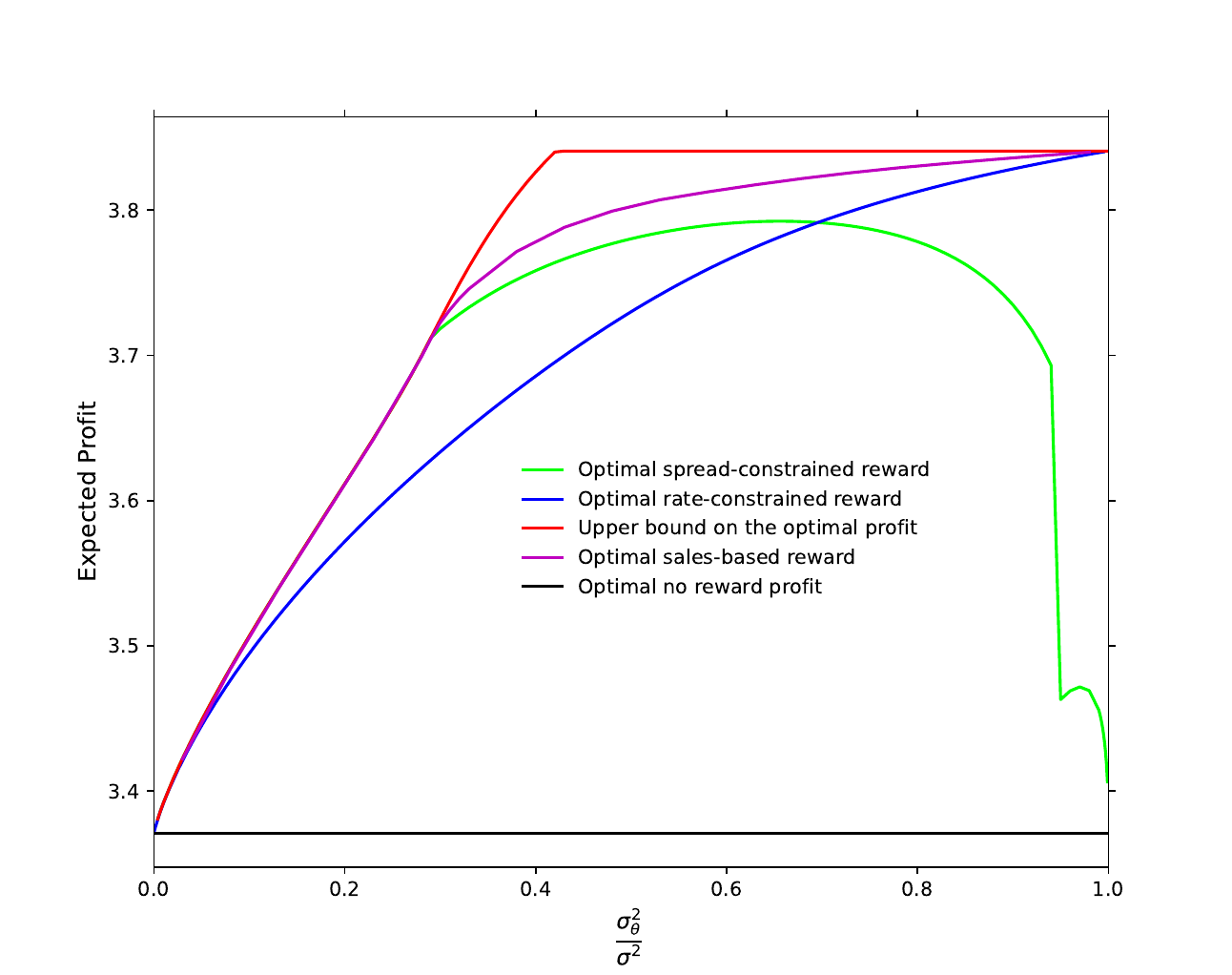}
  \caption{The expected profit for the spread-constrained ($\rh^*_{SC}(\cdot)$) and rate-constrained ($\rh^*_{RC}(\cdot)$) reward programs  characterized in Theorem~\ref{theorem::optimal-BSpread} and \ref{theorem::optimallinear}, versus the upper bound on the optimal expected profit ($\Pi^H$), the optimal expected profit ($\E_v[\Pi(p,\rh^*)]$) from numerical simulations, and the expected profit of the no reward case ($\E_v[\Pi(p,0)]$), for a sample choice of $\theta=5$ and total uncertainty $\sigma$ normalized to 1. The resultant optimal price for the no reward case is $p=3.91$.}
  \label{fig::ProfitsAll}
\end{figure}

We use the upper bound $\Pi^H=\min(\Pi_1^H,\Pi_2^H)$ for the optimal expected profit to evaluate the performance of the reward programs characterized in Theorem~\ref{theorem::optimal-BSpread} and \ref{theorem::optimallinear}.
This upper bound, along with the expected profit resulted from the $\rh^*_{SC}(\cdot)$ and $\rh^*_{RC}(\cdot)$ are plotted in Figure~\ref{fig::ProfitsAll} for a sample choice of $\theta$, and various allocations of the total uncertainty in consumers' valuations $\sigma^2=\sigma_\theta^2+\sigma_\ep^2$ between the two sources of uncertainty (that is, diversity in tastes and uncertainty in common quality).
When most of the uncertainty in valuations comes from the diversity in tastes, the expected profit resulted from  $\rh^*_{SC}(\cdot)$ coincides with the upper bound on the expected profit, certifying the optimality of $\rh^*_{SC}(\cdot)$ and tightness of the upper bound in this regime. On the other hand, as $\sigma_\ep\to 0$, that is when the uncertainty in valuations is primarily due to the uncertainty in the common quality, the expected profit  from  $\rh^*_{RC}(\cdot)$
approaches the upper bound on the expected profit, implying asymptotic optimality of $\rh^*_{RC}(\cdot)$ as $\sigma_\ep\to 0$.

\begin{Theorem}
\label{theorem::performance}
Fix the total uncertainty in consumers' valuations $\sigma^2=\sigma_\theta^2+\sigma_\ep^2$. Then,
\begin{itemize}
\item[i)]\noindent  if $\frac{\sigma_\ep}{\sigma_\theta}\geq\frac{p}{\sqrt{2\pi}\sigma}$, then the spread-constrained reward program $\rh^*_{SC}(\cdot)$ characterized in Theorem~\ref{theorem::optimal-BSpread} is optimal. That is, $\rh^*_{SC}(\cdot)$ is also the solution to the variational optimization problem in \eqref{eq::optimal}.
    %$\E_v[\Pi(p,\rh^*(v))]=\E_v[\Pi(p,\rh^*_{SC}(v))]$
\item[ii)]\noindent the rate-constrained reward program $\rh^*_{RC}(\cdot)$ characterized in Theorem~\ref{theorem::optimallinear} is asymptotically optimal, in the sense that
\begin{align}
\lim_{\sigma_\ep\to0}\frac{\E_v[\Pi(p,\rh^*_{RC}(v))]}{\E_v[\Pi(p,\rh^*(v))]}=1,
\end{align}
where $\rh^*(v)$ is the optimal reward program formulated in \eqref{eq::optimal}.
\end{itemize}
\end{Theorem}
\noindent\textit{Proof.} {See the appendix.}$\hfill\blacksquare$

%two points about the linear one: it is not the limit of $l\to \infty$
%discontinuity of the profit of joint reward: especially that for small diversity in tastes, seller can practically extract the full surplus.

One last point worth mentioning here is that the optimal rate-constrained reward program of Theorem~\ref{theorem::optimallinear}, which is in form of a saturated linear function with the range $[0,p]$, should not be interpreted as the limit of the sequence of optimal reward programs $\{\rh^*_{\sigma_\ep}\}$ as $\sigma_\ep\to 0$. Although Theorem~\ref{theorem::optimal-characteristics} establishes that the number of levels in $\rh^*_{\sigma_\ep}$ goes to infinity as $\sigma_\ep\to 0$, these functions are ``full-refund or nothing'' policies only taking values in $\{0,p\}$. Therefore, the sequence of optimal reward programs $\{\rh^*_{\sigma_\ep}\}$ as $\sigma_\ep\to 0$, do not converge to $\rh^*_{RC}$ and may not even be convergent at all.

\section{Conclusions}
Motivated by the significant gain in profit margins for products and services inherently exhibiting network effects, we ask whether inducing a similar effect via a reward mechanism that conditions the amount of the reward on the sales volume, is beneficial to a firm selling a product with no network effect. 
Consumers in our model are differentiated by their intrinsic valuations, which parametrize both the average common quality of the product and their idiosyncratic taste-related preferences. 
Both the firm and consumers are uncertain about the average quality, sharing a common prior belief on it.   
Each consumer privately observes her own valuation, but cannot separate the common quality from her taste-related component. 
The combination of heterogeneous tastes and uncertainty in the common quality induces a global game with correlated private valuations among consumers.
Consumers' heterogeneous beliefs on sales volume enables  the firm to discriminate net prices over consumers' valuations.
As one of our key results, we show that
while firms can often gain when a product naturally exhibits positive network effect,
creating such an effect using a sales-based reward program could have an adverse effect on their profit. Nevertheless, incentive programs such as group-buying and referral rewards
can be still beneficial due to their effectiveness in fostering word of mouth and social influence, scale economies, and reducing supply-demand mismatch under capacity constraint
in situations discussed in the literature of group-buying and referral reward programs.
Using variational optimization techniques, we identify several key characteristics of the optimal reward program.
The optimal reward program is a ``full-refund or nothing'' policy, fully reimbursing  buyers if and when the realized quality lies in one of the finitely many refund-eligible intervals. As consumers' tastes become less diverse and valuations concentrate further around the quality, the number of these intervals grows unboundedly. Moreover, in this regime,  limiting the number of refund intervals degrades its performance to that of the no-reward case. Motivated by the limitations of the optimal solution, 
we propose two alternatives to the optimal design by analytically solving for the optimal solution within two subspaces of sales-based reward functions:
one with a constraint on the reward spread  and the other with a constraint on its rate of change. 
Despite their simple structures, these two reward programs perform provably-well, asymptotically recovering the optimal solution.

% \bibliographystyle{plainnat}
% \bibliography{refs}

\input{main.bbl}
\bibliographystyle{plainnat}

\newpage
\begin{appendix}
\input{appendix-v7}

\end{appendix}

\end{document}

%% file: header.tex
\usepackage{amsmath, amssymb, amsthm, graphicx} 					% Formerly amstex, advanced math extensions
\usepackage{bm}									
\usepackage{dsfont}		
\usepackage{accents}
\usepackage[usenames,dvipsnames]{color}
\usepackage{setspace}														% allows for a change in line spacing, e.g. double spacing
\usepackage{url}																% \url{...} in which _ & $ % are allowed
\usepackage{natbib}															% allows for Harvard-style citation, e.g. Samuelson (1953) as opposed to Samuleson ([1])
\bibpunct{(}{)}{;}{a}{,}{,}											% costumizes the style of author-year citations, see Get to grips with.... LateXBibliography for details

\usepackage[plainpages=false]{hyperref}					% manage links within the document

\hypersetup{
 colorlinks=true,
 citecolor=Maroon,
 linkcolor=Maroon,
 urlcolor=Maroon}

\usepackage{pst-all}                            % PSTricks package for graphics

\usepackage{nicefrac}														% to have some nicer \frac
\newcommand{\ubar}[1]{\underaccent{\bar}{#1}}

\newcommand{\Prob}{\mbox{Prob}}

\newcommand{\E}{\mathbb{E}}

\newcommand{\ep}{\epsilon}

\usepackage{xifthen}% provides \isempty test

\newcommand{\normcdf}[2]{
\ifthenelse{\isempty{#2}}%
{\Phi(#1)}%
{\Phi(\frac{#1}{#2})}%
}

\newcommand{\normpdf}[2]{
\ifthenelse{\isempty{#2}}%
{\phi(#1)}%
{\frac{\phi(\frac{#1}{#2})}{#2}}%
}

\newtheorem{Theorem}{Theorem}

\newtheorem{Lemma}{Lemma}
\newtheorem*{Proof*}{Proof}

\newtheorem{Assumption}{Assumption}

\def \bi {\begin{itemize}}
\def \ei {\end{itemize}}
\def \be {\begin{equation}}
\def \ee {\end{equation}}
\def \ba {\begin{align}}
\def \ea {\end{align}}

\psset{unit=1mm}

%% file: appendix-v7.tex
% \section{Proofs}
%%%%%%%%%%%%%%%%%%%%%%%%%%%%%%%%%%%%%%%%%%%%%%%%%%%%%%%%%%%%%%%%%%%%%%%%%%%%%
%%%%%%%%%%%%%%%%%%%%%%%%%%%%%%%%%%%%%%%%%%%%%%%%%%%%%%%%%%%%%%%%%%%%%%%%%%%%%
\noindent\textit{Proof of Theorem~\ref{theorem::optimal-characteristics}.} Most of the proof is already given in the body of the paper in Section~\ref{subsec::optimality} and Section~\ref{subsec::performanceA}, so we only fill in the gaps by providing the details where needed.

\noindent\textit{Existence of an optimal reward program:}
Finding the optimal reward program as posed in \eqref{eq::optimal} involves maximizing the expected profit of the seller over $\rh\in L^\infty(\mathbb{R};[0,p])$, that is the closed half-sphere of radius $p$ in $L^\infty(\mathbb{R})$, and $c\in [0,p]$.
$L^\infty(\mathbb{R};[0,p])$ is weak$^*$ compact, according to Alaoglu's theorem (see, e.g., \cite{Optimization_Luenberger}, Page~128).
To prove the existence of a global maximizer for \eqref{eq::optimal}, it thus suffices to show the continuity of $\E_v[(p-\rh(v)) \Phi(\frac{v-c}{\sigma_\ep})]$ in $(c,\rh)$.
Continuity of $\E_v[(p-\rh(v)) \Phi(\frac{v-c}{\sigma_\ep})]$ is apparent from
\begin{align}
|\E_v[(p-\rh_n(v)) \Phi(\frac{v-c_n}{\sigma_\ep})-(p-\rh(v)) \Phi(\frac{v-c}{\sigma_\ep})]|&\leq \frac{\normpdf{0}{}}{\sigma_\ep}|c_n-c|+|\E_v[(\rh_n(v)-\rh(v))\normcdf{v-c}{\sigma_\ep}]|\nonumber\\
&\leq \frac{\normpdf{0}{}}{\sigma_\ep}|c_n-c|+\E_v[|\rh_n(v)-\rh(v)|].
\end{align}

%Compactness of $L_1(\mathbb{R};[0,p],d\normcdf{v-\theta}{\sigma_\theta})$ using the Frechet-Kolmogorov theorem:\\
%%page 34, Ito-long
%i) equi-boundedness:
%\begin{align}
%\sup_{\rh:\mathbb{R}\to[0,p]}\E_v[|\rh(v)|]<p<\infty\nonumber\\
%\end{align}
%ii) equi-continuity:
%\begin{align}
%\lim_{|t|\to 0}\sup_{\rh:\mathbb{R}\to[0,p]}\E_v[|\rh(v+t)-\rh(v)|],
%\end{align}
\noindent\textit{Upper bound on the expected profit for the optimal reward program with $l$ full-refund eligible intervals:}
The inequality in \eqref{eq::refundeligible2} 
provides us with a very useful upper-bound on the mass of qualities in $[\mu_c,\mu_p]$ that are not eligible for a full-refund. We restate this upper bound here:
\begin{align}
\label{eq::refundeligible2-appendix}
|[\mu_c,\mu_p]\setminus\mathcal V_{\rm refund}|<2p\tau\normcdf{-\delta}{\sigma_v}+2l\delta,
\end{align}
for any choice of $\delta>0$. 
We next use this to prove parts i) and ii) of the theorem on the effect of the number of refund-eligible intervals on the profit of the seller.  
The seller fully refunds the buyers if the realized quality $v\in \mathcal V_{\rm refund}$, and charges each buyer a net price equal to $p$ otherwise. Combining this with the above upper bound, we can come up with the following bound on the profit of the seller:
\begin{align}
\label{eq::profuboundL-appendix}
\E_v[(p-\hat r(v)) \Phi(\frac{v-c}{\sigma_\ep})]=&p\int_{\mathbb{R}\setminus\mathcal V_{\rm refund}}\normcdf{v-c}{\sigma_\ep}\normpdf{\theta-v}{\sigma_\theta}dv\nonumber\\
<&p\int_{[\mu_c,\mu_p]\setminus\mathcal V_{\rm refund}}\normcdf{v-c}{\sigma_\ep}\normpdf{\theta-v}{\sigma_\theta}dv
+ p\int_{\mathbb{R}\setminus[\mu_c,\mu_p]}\normcdf{v-c}{\sigma_\ep}\normpdf{\theta-v}{\sigma_\theta}dv\\
<&\frac{p}{\sqrt{2\pi}\sigma_\theta}(2p\tau\normcdf{-\delta}{\sigma_v}+2l\delta)
+ p\int_{\mathbb{R}\setminus[\mu_c,\mu_p]}\normcdf{v-c}{\sigma_\ep}\normpdf{\theta-v}{\sigma_\theta}dv
\end{align}
Fix the total uncertainty in consumers' valuations $\sigma^2=\sigma_\theta^2+\sigma_\ep^2$.
When $\sigma_\epsilon\to0$ the second term in the above upper bound will have an asymptotic value of $p\normcdf{\theta-p}{\sigma}$. To see this, notice that
\begin{align}
\label{eq::profoutbound}
\int_{\mathbb{R}\setminus[\mu_c,\mu_p]}\normcdf{v-c}{\sigma_\ep}\normpdf{\theta-v}{\sigma_\theta}dv&<
\normcdf{\theta-\mu_p}{\sigma_\theta}+\frac{1}{\sqrt{2\pi}\sigma_\theta}\int_{-\infty}^{\mu_c}\normcdf{v-c}{\sigma_\ep}dv\nonumber\\
&<\normcdf{\theta-\mu_p}{\sigma_\theta}+\frac{\sigma_\ep}{\sqrt{2\pi}\sigma_\theta}(\sigma_\ep(\theta-c)\normcdf{\sigma_\ep(\theta-c)}{}+\normpdf{\sigma_\ep(\theta-c)}{}),
\end{align}
which is asymptotically upper-bounded by $\normcdf{\theta-p}{\sigma}$ as $\sigma_\ep\to0$. Using this along with \eqref{eq::profuboundL-appendix}, we then get:
\begin{align}
\label{eq::uboundL2-appendix}
\liminf_{\sigma_\epsilon\to0}\E_v[\Pi(p,\hat r(v))]&\leq p\normcdf{\theta-p}{\sigma}+
\liminf_{\sigma_\epsilon\to0}{\frac{p}{\sqrt{2\pi}\sigma_\theta}(2p\tau\normcdf{-\delta}{\sigma_v}+2l\delta)}\nonumber\\
&\leq p\normcdf{\theta-p}{\sigma}+\liminf_{\sigma_\epsilon\to0}{\frac{p}{\sqrt{2\pi}\sigma}(2p\normcdf{-\delta}{\sigma_\epsilon}+2l\delta)}.
\end{align}
Let us choose $\delta = \sqrt{\frac{\sigma_\epsilon}{l}}$. Then, 
\begin{align}
\liminf_{\sigma_\epsilon\to0}\E_v[\Pi(p,\hat r(v))]&\leq p\normcdf{\theta-p}{\sigma}+
\liminf_{\sigma_\epsilon\to0}{\frac{p}{\sqrt{2\pi}\sigma}(2p\normcdf{-1}{\sqrt{l\sigma_\epsilon}}+2\sqrt{l\sigma_\epsilon })}.
\end{align}
For fixed $l$, this implies that
\begin{align}
\liminf_{\sigma_\epsilon\to0}\E_v[\Pi(p,\hat r(v))]&\leq p\normcdf{\theta-p}{\sigma},
\end{align}
where  $p\normcdf{\theta-p}{\sigma}$ is the expected profit of the no reward case. This completes the proof of part ii) of the theorem.
On the other hand, we can find the exact value of the asymptotic optimal profit:
\begin{align}
\E_v[\Pi(p,\hat r^*(v))]\to\int_{0}^{p}x\normpdf{\theta-x}{\sigma}dx+p\normcdf{\theta-p}{\sigma}\quad \text{as} \sigma_\epsilon\to 0.
\end{align}
This is because, as we will show later in 
Theorem~\ref{theorem::performance}, the expected profit from $\rh^*_{RC}(\cdot)$ approaches $\Pi_1^H$ as $\sigma_\ep\to0$ where the upper bound $\Pi_1^H$ given by \eqref{eq::profubound1} has the above value. Using \eqref{eq::uboundL2-appendix} with the optimal solution $\rh^*$ as $\rh$, it should thus hold that
\begin{align}
\label{eq::performance3}
\liminf_{\sigma_\epsilon\to0}{\frac{p}{\sqrt{2\pi}\sigma}(2p\normcdf{-\delta}{\sigma_\epsilon}+2l\delta)}\geq\int_{0}^{p}x\normpdf{\theta-x}{\sigma}dx.
\end{align}
Choosing the same $\delta = \sqrt{\frac{\sigma_\epsilon}{l}}$, we get $\liminf_{\sigma_\epsilon\to0} \sigma_\epsilon l>0$, meaning $\frac{1}{\sigma_\epsilon}=\mathcal O(l)$, or equivalently $l=\Omega(\frac{1}{\sigma_\epsilon})$, completing the proof.
We can even push further and optimize over $\delta$ by minimizing $2p\normcdf{-\delta}{\sigma_\epsilon}+2l\delta$. Firs order condition gives $\phi\big(\frac{\delta}{\sigma_\epsilon}\big)= \frac{l\sigma_\epsilon}{p}$. This has a solution if and only if $l\sigma_\epsilon\leq \frac{p}{\sqrt{2\pi}}$, in which case we choose it (the positive one) as $\delta$, and we choose $\delta=0$ otherwise.
When $l\sigma_\epsilon\leq \frac{p}{\sqrt{2\pi}}$, we get
\begin{align}
2p\normcdf{-\delta}{\sigma_\epsilon}+2l\delta = 2p\Big(\normcdf{-\delta}{\sigma_\epsilon}+\frac{\delta}{\sigma_\epsilon}\phi\big(\frac{\delta}{\sigma_\epsilon}\big)\Big).
\end{align}
Let $\kappa$ be the unique solution to 
\begin{align}
\normcdf{-\kappa}{}+\kappa\normpdf{\kappa}{} = \frac{\sqrt{2\pi}\sigma}{2p^2}\int_{0}^{p}x\normpdf{\theta-x}{\sigma}dx.
\end{align}
Using \eqref{eq::performance3} and that in the case where $l\sigma_\epsilon>\frac{p}{\sqrt{2\pi}}$ we already have $l\sigma_\epsilon>{p}\normpdf{0}{}$,
we reach at
\begin{align}
\liminf_{\sigma_\epsilon\to0} l \sigma_\epsilon\geq p\normpdf{\kappa}{}.
\end{align}
~$\hfill\blacksquare$

\noindent\textit{Proof of Lemma~\ref{lemma::necessarycond}.}
Proof follows directly from \eqref{eq::EPi}.
$\hfill\blacksquare$

%%%%%%%%%%%%%%%%%%%%%%%%%%%%%%%%%%%%%%%%%%%%%%%%%%%%%%%%%%%%%%%%%%%%%%%%%
%\noindent\textit{Proof of Lemma~\ref{lemma::monotoneBNE}.}
%Proof is given in the body of the paper right above the lemma.
%$\hfill\blacksquare$

%%%%%%%%%%%%%%%%%%%%%%%%%%%%%%%%%%%%%%%%%%%%%%%%%%%%%%%%%%%%%%%%%%%%%%%%%%%%%
%To prove \eqref{eq::monotonereward} for a decreasing reward function $r(\cdot)$, it suffices to show that the posterior beliefs satisfy the first order stochastic dominance property, that is
%$P_{v|v_1}[v\geq x]\geq P_{v|v_2}[v\geq x]$ for every $x\in\mathbb{R}$ and $v_1\geq v_2$.
%This then follows from the
%monotone likelihood ratio property (MLRP) of the normal distribution (see \cite{Milgrom_81}). The proof for the other case is similar.

%%%%%%%%%%%%%%%%%%%%%%%%%%%%%%%%%%%%%%%%%%%%%%%%%%%%%%%%%%%%%%%%%%%%%%%%%%%%%

\noindent\textit{Proof of Theorem~\ref{theorem::monotonereward}.}
We first prove the following lemma on the monotonicity of expected rewards for monotone reward programs.
\begin{Lemma}
\label{lemma::monotonereward}
Suppose the consumers follow a threshold strategy of the form $a_i={\bf 1}\{v_i>c\}$. Then, for any increasing non-constant reward program $r(\cdot)$ the expected rewards of the consumers are strictly increasing with their valuations. In other words,
\begin{equation}
\label{eq::monotonereward}
v_1> v_2\Rightarrow\E_{v|v_1}[\rh(v)]> \E_{v|v_2}[\rh(v)].
\end{equation}
Similarly, the expected rewards are decreasing with consumers' valuations for any decreasing reward program.
\end{Lemma}
%\noindent\textit{Proof.} {See the appendix.}~$\hfill\blacksquare$
\noindent\textit{Proof.}
Recall that $v|x\sim N(\tau x+(1-\tau)\theta,\sigma_v^2)$, where $\tau=\frac{\sigma_\theta^2}{\sigma_\ep^2+\sigma_\theta^2}$ and $\sigma_v^2=\frac{\sigma_\ep^2\sigma_\theta^2}{\sigma_\ep^2+\sigma_\theta^2}$. Therefore,
\begin{align}
\E_{v|x}[\rh(v)]&=\int_{-\infty}^{\infty}{\rh(v)\frac{\phi(\frac{v-\tau x-(1-\tau)\theta}{\sigma_v})}{\sigma_v}dv}\nonumber\\
&=\int_{-\infty}^{\infty}{\rh(v+\tau x)\frac{\phi(\frac{v-(1-\tau)\theta}{\sigma_v})}{\sigma_v}dv},
\end{align}
from which the lemma immediately follows.%~$\hfill\blacksquare$

Based on the above lemma, any increasing non-constant reward program violates condition \eqref{eq::necessarycond} in Lemma~\ref{lemma::necessarycond}, and hence cannot be profitable.
~$\hfill\blacksquare$

%%%%%%%%%%%%%%%%%%%%%%%%%%%%%%%%%%%%%%%%%%%%%%%%%%%%%%%%%%%%%%%%%%%%%%%%%%%%%%%%%%%%%%%%%%%%%%%%%%%%%
\noindent\textit{Proof of Lemma~\ref{lemma::rewardspread}.}
For $\frac{\partial}{\partial v_i}\E_v[u_i|v_i]\geq0$ to hold, it suffices to have $|\frac{\partial}{\partial v_i}\E_{v|v_i}[\rh(v)]|\leq1$. Let $\mu_i=\tau v_i+(1-\tau)\theta$. Then,
\begin{align}
\frac{\partial}{\partial v_i}\E_{v|v_i}[\rh(v)]&=
\E_{v|v_i}[\rh(v)\frac{-(v-\mu_i)}{\sigma_v}\frac{\partial}{\partial v_i}(\frac{v-\mu_i}{\sigma_v})]\nonumber\\
&=\frac{\tau}{\sigma_v}\E_{v|v_i}[\rh(v)\frac{(v-\mu_i)}{\sigma_v}].
\end{align}
Using this, it is easy to see that
\begin{align}
|\frac{\partial}{\partial v_i}\E_{v|v_i}[\rh(v)]|\leq
\frac{\tau}{\sigma_v}(r_{\max}-r_{\min})\int_{0}^{\infty}{x\phi(x)dx}=\frac{\tau}{\sqrt{2\pi}\sigma_v}(r_{\max}-r_{\min}).
\end{align}
Therefore, if $r_{\max}-r_{\min}\leq \frac{\sqrt{2\pi}\sigma_v}{\tau}=\sqrt{2\pi}\frac{\sigma_\ep\sigma}{\sigma_\theta}$ then $\frac{\partial}{\partial v_i}\E_v[u_i|v_i]\geq0$, which completes the proof.$\hfill\blacksquare$

\noindent\textit{Proof of Theorem~\ref{theorem::optimal-BSpread}.}
WLOG, we normalize the total uncertainty in consumers' valuations to 1, that is, we assume $\sigma^2=\sigma_\ep^2+\sigma_\theta^2=1$.
We prove the theorem in a few steps: We first fix a window $[r_{\min}, r_{\max}]$ of length $r_M$ (i.e., $r_{\max}-r_{\min}=r_M$) and maximize the expected profit over $\hat r:\mathbb{R}\to[r_{\min},r_{\max}]$ and $c\in\mathbb{R}$. We then optimize over $r_{\min}$, where we show that the expected profit is indeed maximized when $r_{\min}=0$. The latter is done using Lemma~\ref{lemma::logconcave1}.

Fix a number $r_{\min}\in [0,p-r_M]$, and let $r_{\max}=r_{\min}+r_M$. 
To maximize the expected profit for $\hat r:\mathbb{R}\to[r_{\min},r_{\max}]$ and $c\in\mathbb{R}$, we use the Lagrangian 
\begin{align}
L(\rh,c,\lambda)=\E_v[\Pi(p,\rh(v))]-\lambda (p-c-\E_{v|v_i=c}[\hat r(v)]),
\end{align}
and find the pair $(\hat r,c)$ where $c\in\mathbb{R}$ and $\hat r:\mathbb{R}\to[r_{\min},r_{\max}]$ maximizing the Lagrangian.
Recalling that $\E_v[\Pi(p,\rh(v))]=\int_{c}^{\infty}(p-\E_{v|v_i}[\hat r(v)])\psi(v_i)dv_i$, we can write the
first order condition for optimal cutoff $c$ as
\begin{align}
0=-(p-\E_{v|v_i=c}[\hat r(v)])\psi(c)-\lambda (-1+\frac{\tau}{\sigma_v}\E_{v|v_i=c}[(\frac{\mu_c-v}{\sigma_v})\hat r(v)]),
\end{align}
which yields
\begin{align}
\label{eq::FOC1B}
\lambda=\frac{c\phi(\theta-c)}{1-\frac{\tau}{\sigma_v}\E_{v|v_i=c}[(\frac{\mu_c-v}{\sigma_v})\hat r(v)]}.
\end{align}

To find the optimal value of $\hat r(v)$ at each realization $v$ of the quality,  we use $\Pi(p,\rh(v))=(p-\rh(v))\normcdf{v-c}{\sigma_\ep}$ to
find the weight of $\rh(v)$ in the Lagrangian $L$:
\be
\label{eq::rv-weight}
-\normcdf{v-c}{\sigma_\ep}\normpdf{v-\theta}{\sigma_\theta}+\lambda \normpdf{v-\mu_c}{\sigma_v}.
\ee
Since $\hat r:\mathbb{R}\to[r_{\min},r_{\max}]$, the reward $\hat r(v)$ maximizing $L$ is of the form
\be
\label{eq:rminmax}
\hat r(v)=
\begin{cases}
r_{\min}, & -\normcdf{v-c}{\sigma_\ep}\normpdf{v-\theta}{\sigma_\theta}+\lambda \normpdf{v-\mu_c}{\sigma_v}\leq0\\
r_{\max}, & -\normcdf{v-c}{\sigma_\ep}\normpdf{v-\theta}{\sigma_\theta}+\lambda \normpdf{v-\mu_c}{\sigma_v}>0.
\end{cases}
\ee
Let $v_c$ be a zero of \eqref{eq::rv-weight}. Using the identity
\begin{align}
\label{eq::bayes}
\frac{\frac{\phi(\frac{v_c-\mu_c}{\sigma_v})}{\sigma_v}}{\frac{\phi(\frac{v_c-\theta}{\sigma_\theta})}{\sigma_\theta}}=\frac{\frac{\phi(\frac{v_c-c}{\sigma_\ep})}{\sigma_\ep}}{{\frac{\phi(\frac{\theta-c}{\sqrt{\sigma_\ep^2+\sigma_\theta^2}})}{\sqrt{\sigma_\ep^2+\sigma_\theta^2}}}},
\end{align}
and the assumption $\sigma_\ep^2+\sigma_\theta^2=1$, we can reach at
\be
\label{eq::vc}
\frac{\lambda}{\normpdf{\theta-c}{}}=\frac{\normcdf{v_c-c}{\sigma_\ep}}{\normpdf{v_c-c}{\sigma_\ep}},
\ee
showing the uniqueness of $v_c$ given the log-concavity of the normal distribution. The uniqueness of $v_c$ implies that the weight in \eqref{eq::rv-weight} is negative for $v>v_c$ and is positive for $v<v_c$. Putting this together with \eqref{eq:rminmax}, it then follows that for the optimal $\hat r:\mathbb{R}\to[r_{\min},r_{\max}]$,
we have $\rh(v)=r_{\max}$ for $v<v_c$ and $\rh(v)=r_{\min}$ for $v\geq v_c$.
Noting that $\E_{v|v_i=c}[\frac{\mu_c-v}{\sigma_v}]=0$, first order optimality condition for $c$ given in \eqref{eq::FOC1B} is insensitive to $r_{\min}$, and combined with \eqref{eq::vc} leads to
\be
\frac{\lambda}{\normpdf{\theta-c}{}}=\frac{\normcdf{v_c-c}{\sigma_\ep}}{\normpdf{v_c-c}{\sigma_\ep}}=\frac{c}{1-\frac{\tau}{\sigma_v}r_M\normpdf{\frac{v_c-\mu_c}{\sigma_v}}{}}.
\ee
The reward $r_{\min}$ works as a discount as it is paid at all values of $v$. Note that although no discount is profitable when there is no reward ($p$ is the optimal price with no reward), this does not necessarily imply $r_{\min}=0$ when we allow for a nonzero reward program. Using \eqref{eq::rv-weight} the total weight of $r_{\min}$ in $L$ is $-\normcdf{\theta-c}{}+\lambda$. Therefore,
\be
\label{eq:rmin}
r_{\min}=
\begin{cases}
0, & \lambda<\Phi(\theta-c)\\
p-r_M, & \lambda>\Phi(\theta-c)\\
\in[0,p-r_M], & \lambda=\Phi(\theta-c).
\end{cases}
\ee

We claim, however, that having a discount is still not profitable when a nonzero reward is allowed.
For the case $p\leq\sqrt{2\pi}\frac{\sigma_\ep\sigma}{\sigma_\theta}$, we have $r_M=\min(p,\sqrt{2\pi}\frac{\sigma_\ep\sigma}{\sigma_\theta})=p$ yielding $r_{\min}=0$ from \eqref{eq:rmin}. Therefore, we only need to deal with the case where $r_M=\sqrt{2\pi}\frac{\sigma_\ep\sigma}{\sigma_\theta}<p$. We prove $r_{\min}=0$
 by showing that $\lambda<\Phi(\theta-c)$. Suppose (by contradiction) that $\lambda\geq\Phi(\theta-c)$. Then,
\be
\label{eq::FOCboth}
\frac{\Phi(\theta-c)}{\phi(\theta-c)}\leq\frac{\lambda}{\normpdf{\theta-c}{}}=\frac{\normcdf{v_c-c}{\sigma_\ep}}{\normpdf{v_c-c}{\sigma_\ep}}=\frac{c}{1-\sqrt{2\pi}\phi(\frac{v_c-\mu_c}{\sigma_v})}.
\ee
This clearly requires $\frac{v_c-c}{\sigma_\ep}\geq\theta-c$. This, combined with the identity
\be
\label{eq::identity1}
\frac{v_c-c}{\sigma_\ep}=\sigma_\ep(\theta-c)+\sigma_\theta(\frac{v_c-\mu_c}{\sigma_v}),
\ee
implies that
\be
\label{eq::xyz1}
\theta-c\leq\frac{v_c-c}{\sigma_\ep}\leq\frac{1+\sigma_\ep}{\sigma_\theta}(\frac{v_c-\mu_c}{\sigma_v}).
\ee
On the other hand, using $\frac{\Phi(\theta-c)}{\phi(\theta-c)}\geq\frac{\Phi(\theta-p)}{\phi(\theta-p)}=p$ which follows from the optimality of $p$, we can find
\be
\label{eq::FOC2}
\frac{\lambda}{\normpdf{\theta-c}{}}=\frac{c}{1-\sqrt{2\pi}\phi(\frac{v_c-\mu_c}{\sigma_v})}\geq c+r_M\normcdf{v_c-\mu_c}{\sigma_v}\geq\frac{r_M\normcdf{v_c-\mu_c}{\sigma_v}}{\sqrt{2\pi}\phi(\frac{v_c-\mu_c}{\sigma_v})}.
\ee
Comparing \eqref{eq::FOCboth} and \eqref{eq::FOC2}, we can observe that $\frac{\normcdf{v_c-c}{\sigma_\ep}}{\normpdf{\frac{v_c-c}{\sigma_\ep}}{}}\geq\frac{\normcdf{v_c-\mu_c}{\sigma_v}}{\normpdf{\frac{v_c-\mu_c}{\sigma_v}}{}}$,
yielding $\frac{v_c-c}{\sigma_\ep}\geq\frac{v_c-\mu_c}{\sigma_v}$. Incorporating in \eqref{eq::identity1}, we get
\be
\label{eq::xyz2}
\frac{v_c-\mu_c}{\sigma_v}\leq\frac{v_c-c}{\sigma_\ep}\leq\frac{1+\sigma_\theta}{\sigma_\ep}(\theta-c).
\ee
One immediate consequence of \eqref{eq::xyz1} and \eqref{eq::xyz2} is that $\theta-c, \frac{v_c-\mu_c}{\sigma_v},\frac{v_c-c}{\sigma_\ep}\geq0$.

We next present the following useful lemma.
\begin{Lemma}
\label{lemma::logconcave1}
Let $0\leq a\leq b\leq c\leq d$ be such that $a^2+d^2=b^2+c^2$. Then,
\be
\frac{\Phi(a)}{\phi(a)}\times \frac{\Phi(d)}{\phi(d)}\leq \frac{\Phi(b)}{\phi(b)}\times \frac{\Phi(c)}{\phi(c)}.
\ee
\end{Lemma}
\noindent\textit{ Proof}. This easily follows from the log-concavity of $\Phi(\sqrt{x})$ and that ${\phi(a)}\phi(d)=\phi(b)\phi(c)$.%$\hfill\blacksquare$

Applying the above lemma, we can reach at
\be
\frac{\normcdf{v_c-c}{\sigma_\ep}}{\normpdf{v_c-c}{\sigma_\ep}}\times\frac{\normcdf{|v_c-\theta|}{\sigma_\theta}}{\normpdf{|v_c-\theta|}{\sigma_\theta}}\leq
\frac{\normcdf{\theta-c}{}}{\normpdf{\theta-c}{}}\times\frac{\normcdf{v_c-\mu_c}{\sigma_v}}{\normpdf{v_c-\mu_c}{\sigma_v}}.
\ee
From \eqref{eq::FOCboth}, $\frac{\normcdf{v_c-c}{\sigma_\ep}}{\normpdf{v_c-c}{\sigma_\ep}}\geq \frac{\normcdf{\theta-c}{}}{\normpdf{\theta-c}{}}$, and hence,
\be
\frac{\normcdf{v_c-\mu_c}{\sigma_v}}{\normpdf{v_c-\mu_c}{\sigma_v}}\geq\frac{\normcdf{|v_c-\theta|}{\sigma_\theta}}{\normpdf{|v_c-\theta|}{\sigma_\theta}}\geq\sqrt{\frac{\pi}{2}}\sigma_\theta.
\ee
Therefore,
\be
\label{eq::sep}
\frac{\normcdf{v_c-\mu_c}{\sigma_v}}{\normpdf{\frac{v_c-\mu_c}{\sigma_v}}{}}\geq\frac{1}{\sigma_\ep}\sqrt{\frac{\pi}{2}}.
\ee
Using this in \eqref{eq::FOC2}, we get
\be
p\geq c+r_M\normcdf{v_c-\mu_c}{\sigma_v}\geq\frac{1}{\sigma_\theta}\sqrt{\frac{\pi}{2}}\geq\sqrt{\frac{\pi}{2}},
\ee
putting $p$ and hence $\theta$ in the regime $\theta\geq p\geq\sqrt{\frac{\pi}{2}}$. Another useful inequality can be obtained using the log-convexity of the function $h(x)=\frac{\Phi(\theta-x)}{x\phi(\theta-x)}$ for $x\in [0, p]$:
\be
\label{eq::thetac}
\frac{\Phi(\theta-c)}{c\phi(\theta-c)}\geq e^{(p-c)(\frac{2}{p}+\theta-p)}\geq e^{r_M\normcdf{v_c-\mu_c}{\sigma_v}(\frac{2}{p}+\theta-p)}.
\ee

Let $y:=\frac{v_c-\mu_c}{\sigma_v}$. Incorporating $\theta-c\geq p-c\geq r_M\Phi(y)$ in \eqref{eq::xyz1}, we can obtain
\be
{\sigma_\ep}\leq\frac{(1+\sigma_\ep)y}{\sqrt{2\pi}\Phi(y)}\leq \frac{2y}{\sqrt{2\pi}\Phi(y)},
\ee
which along with \eqref{eq::sep} leads to
\be
\frac{y}{\phi(y)}\geq\frac{\pi}{2},
\ee
requiring $y>\frac{1}{2}$. Finally, \eqref{eq::thetac} and \eqref{eq::FOCboth}, along with the fact that $\frac{2}{p}+\theta-p>\frac{3}{2}$ for $\theta\geq\sqrt{\frac{\pi}{2}}$ and optimal price $p=\frac{\Phi(\theta-p)}{\phi(\theta-p)}$, we can show
\be
\frac{1}{1-\sqrt{2\pi}\phi(y)}\geq e^{3\sqrt{\frac{\pi}{2}}\frac{\sigma_\ep}{\sqrt{1-\sigma_\ep^2}}\Phi(y)}.
\ee
With a bit of manipulation, we can verify that the above cannot hold for $y>\frac{1}{2}$ and $\sigma_\ep\geq\sqrt{\frac{\pi}{2}}\frac{\phi(y)}{\Phi(y)}$. This completes the proof of $\lambda<\Phi(\theta-c)$, implying that $r_{\min}=0$.

The last point to be made is about the uniqueness of the equilibrium threshold strategy for the subgame among the consumers under the optimal spread-constrained reward program. This is important, as in case of multiple equilibria seller needs to use other measures to speculate which equilibrium is more likely to be followed by consumers. More generally, we can show the uniqueness of the cutoff $c$ for any decreasing reward program $r(\normcdf{v-c}{\sigma_\ep})$. To prove, we write the indifference equation as
\begin{align}
\label{eq::theorem1::tmp1}
p=c+\E_{v|v_i=c}[r(\Phi(\frac{v-c}{\sigma_\ep}))]=c+\E_{v|v_i=p}[r(\Phi(\frac{v-(1-\tau)c-\tau p}{\sigma_\ep}))],
\end{align}
where the second equality follows from combining  $v|v_i=c\sim N(\tau c+(1-\tau)\theta,\sigma_v^2)$ and $v|v_i=p\sim N(\tau p+(1-\tau)\theta,\sigma_v^2)$. It is easy to see  that the RHS in \eqref{eq::theorem1::tmp1} is strictly increasing in $c$, implying the uniqueness of $c$.
~$\hfill\blacksquare$

%%%%%%%%%%%%%%%%%%%%%%%%%%%%%%%%%%%%%%%%%%%%%%%%%%%%%%%%%%%%%%%%%%%%%%%%%%%%%
\noindent\textit{Proof of Lemma~\ref{lemma::rewardROC}.}
Similar to the proof of Lemma~\ref{lemma::rewardspread}, it suffices to show that $|\frac{\partial}{\partial v_i}\E_{v|v_i}[\rh(v)]|\leq1$. This immediately follows from
\begin{align}
\E_{v|v_i}[\rh(v)]=\int_{\mathbb{R}}{\rh(v)\normpdf{v-\mu_i}{\sigma_v}dv}=\int_{\mathbb{R}}{\rh(v+\mu_i)\normpdf{v}{\sigma_v}dv},
\end{align}
noting that $|\frac{\partial}{\partial v_i}\rh(v+\mu_i)|\leq \frac{1}{\tau}\frac{\partial \mu_i}{\partial v_i}=1$.
$\hfill\blacksquare$

%%%%%%%%%%%%%%%%%%%%%%%%%%%%%%%%%%%%%%%%%%%%%%%%%%%%%%%%%%%%%%%%%%%%%%%%%%%%%

\noindent\textit{Proof of Theorem~\ref{theorem::optimallinear}.} Without loss of generality, we normalize the total uncertainty in consumers' valuations to 1, that is, we assume $\sigma^2=\sigma_\ep^2+\sigma_\theta^2=1$.
Denote with $PC(\mathbb{R};[a,b])$ the space of piecewise continuous functions from $\mathbb{R}$ to $[a,b]$.
Then, $\frac{d}{dv}\rh(v)=-u(v)$ for some $u\in PC(\mathbb{R};[-\frac{1}{\tau},\frac{1}{\tau}])$, and
\begin{equation}
\label{eq::rh2u}
\rh(v)=r_0-\int_{-\infty}^{v} u(w) dw,
\end{equation}
where $r_0:=r(\bar a(0))=\rh(-\infty)$.
We can reformulate the problem of finding the optimal reward program with the rate constraint as
\begin{gather}
\label{eq::linearopt}
\underset{\substack{u\in PC(\mathbb{R};[-\frac{1}{\tau},\frac{1}{\tau}]), c,r_0\in\mathbb{R}}}{\text{maximize }} \E_v[\Pi(p,\rh(v))],\\
\text{subject to:}\\
0\leq \hat r(v)\leq p, \text{ for all } v\in\mathbb{R},\\
c+\E_{v|v_i=c}[\rh(v)]=p.
\end{gather}

We proceed by first establishing the existence of Lagrange multipliers for the above optimization problem. 
We then recast the Lagrangian as a function of $u$ instead of $\hat r$ using integration by part and  that $\frac{d}{dv}\rh(v)=-u(v)$. 
We then use the optimality conditions and the complementarity slackness property to establish that while $\hat r^*$ is off the boundary, then $u^*=\frac{1}{\tau}$. This enables us to parameterize the optimal solution by a triplet $(c, v_L, v_H)\in \mathbb{R}^3$, where
the optimal reward stays at $\hat r^*(v)=p$ for $v\leq v_L$, and at $\hat r^*(v)=0$ for $v\geq v_H$, and decreases at the fixed rate $\frac{1}{\tau}$ for $v\in [v_L, v_H]$.
We then solve for $(c, v_L, v_H)$ by putting together the optimality condition, the indifference equation, and that $v_H-v_L = p\tau$.

The existence of Lagrange multipliers for this problem can be established using a regularity condition, which basically requires the linearized approximation of the constraint space around the optimal solution to have a feasible interior point (see, e.g., Definition 1.5 in \cite{Ito_2008} for the explicit statement of the regularity used here).
Denote the optimal cutoff associated with \eqref{eq::linearopt} by $c^*$. One can easily verify that for $\rh=p-c^*$ and $c=c^*$ all the inequalities are strict and the indifference equation is satisfied. Since the only nonlinearity in \eqref{eq::linearopt} (that is, cutoff) is kept unchanged, this verifies the regularity condition discussed above.
Let the Lagrangian be
\begin{align}
\label{eq::L1}
L(u,r_0,c,\lambda,\mu_1,\mu_2)=&\E_v[\Pi(p,\rh(v))]+\lambda (c+\E_{v|v_i=c}[\hat r(v)]-p)\nonumber\\
&+\int_{\mathbb{R}}{\hat r(v)d\mu_1(v)}+\int_{\mathbb{R}}{(p-\hat r(v))d\mu_2(v)},
\end{align}
where $\lambda\in\mathbb{R}$, and $\mu_1,\mu_2\in BV(\mathbb{R})$ are increasing upper-semicontinuous functions with bounded variation vanishing at $+\infty$. Denote the optimal solution of \eqref{eq::linearopt} with $(u^*,r_0^*,c^*)$ and the resulting optimal reward with $\rh^*(\cdot)$. Optimality conditions given by Lagrange theorem are then
\begin{gather}
\label{eq::FOC1a}
\frac{\partial}{\partial c}L(u^*,r_0^*,c,\lambda,\mu_1,\mu_2)_{|c=c^*}=0,\\
L(u^*,r_0^*,c^*,\lambda,\mu_1,\mu_2)\geq L(u,r_0^*,c^*,\lambda,\mu_1,\mu_2),\text{ for all }u\in PC(\mathbb{R};[-\frac{1}{\tau},\frac{1}{\tau}]),\\
L(u^*,r_0^*,c^*,\lambda,\mu_1,\mu_2)\geq L(u^*,r_0,c^*,\lambda,\mu_1,\mu_2),\text{ for all }r_0\in\mathbb{R},
\end{gather}
and the complementary slackness property requires
\begin{align}
\label{eq::CS1}
\int_{\mathbb{R}}{\hat r^*(v)d\mu_1(v)}=&0,\nonumber\\
\int_{\mathbb{R}}{(p-\hat r^*(v))d\mu_2(v)}=&0,
\end{align}
implying that $\mu_1(v)$ and $\mu_2(v)$ can only change when $\rh$ is on the boundary.

Integration by part for each term in RHS of \eqref{eq::L1} using $\frac{d}{dv}\rh(v)=-u(v)$, we can recast the Lagrangian as
\begin{align}
\label{eq::L2}
L(u,r_0,c,\lambda,\mu_1,\mu_2)=&(p-r_0)\Phi(\theta-c)+\int_{\mathbb{R}}u(v)\int_{v}^{+\infty}{\normcdf{w-c}{\sigma_\ep}\normpdf{\theta-w}{\sigma_\theta}dw}dv\\
&+\lambda (c-p+r_0-\int_{\mathbb{R}}{u(v)\normcdf{\mu_c-v}{\sigma_v}dv})\nonumber\\
&+\int_{\mathbb{R}}{u(v)\mu_1(v) dv}-\int_{\mathbb{R}}{u(v)\mu_2(v)dv}-r_0\mu_1^0-(p-r_0)\mu_2^0,
\end{align}
where $\mu_1^0=\mu_1(-\infty)$ and $\mu_2^0=\mu_2(-\infty)$. Regrouping the terms, we get
\begin{align}
\label{eq::L3}
L(u,r_0,c,\lambda,\mu_1,\mu_2)=&(p-r_0)\Phi(\theta-c)+\lambda (c-p+r_0)-r_0\mu_1^0-(p-r_0)\mu_2^0\nonumber\\
&+\int_{\mathbb{R}}u(v)(G(v)+\mu_1(v)-\mu_2(v))dv,
\end{align}
where
\begin{align}
\label{eq::G(v)}
G(v)=\int_{v}^{+\infty}{\normcdf{w-c}{\sigma_\ep}\normpdf{\theta-w}{\sigma_\theta}dw}-\lambda\normcdf{\mu_c-v}{\sigma_v}.
\end{align}
It is easy to see that
\begin{align}
\frac{\frac{\partial}{\partial v}G(v)}{\normpdf{\mu_c-v}{\sigma_v}}=\lambda-\frac{\normcdf{v-c}{\sigma_\ep}}{\normpdf{v-c}{\sigma_\ep}}\phi(\theta-c),
\end{align}
which is a decreasing function of $v$, thus having a unique root $v_c$ for $\lambda>0$ and no root otherwise.
This implies that $G(v)$ is increasing for $v<v_c$ and decreasing for $v\geq v_c$, where $v_c$ is the unique solution of
\begin{align}
\frac{\lambda}{\phi(\theta-c)}=\frac{\normcdf{v-c}{\sigma_\ep}}{\normpdf{v-c}{\sigma_\ep}},
\end{align}
for $\lambda>0$, and $v_c=-\infty$ if $\lambda\leq0$.

From the optimality condition  \eqref{eq::FOC1a} and complementary slackness property \eqref{eq::CS1}, we can obtain
\begin{align}
\label{eq::uoptimal}
u^*(v)=\begin{cases}
\frac{1}{\tau}, &G(v)+\mu_1(v)-\mu_2(v)>0\\
-\frac{1}{\tau}, &G(v)+\mu_1(v)-\mu_2(v)<0\\
0, &G(v)+\mu_1(v)-\mu_2(v)=0.
\end{cases}
\end{align}
Next, we claim that $G(v)+\mu_1(v)-\mu_2(v)\geq0$ for all $v\in\mathbb{R}$.
Suppose, by contradiction, that $G(v_0)+\mu_1(v_0)-\mu_2(v_0)<0$ for some $v_0\geq v_c$ (the case $v_0<v_c$ is similar).
Using the fact that $\mu_2(v)$ is increasing and $\mu_1(v)$ is fixed while off the boundary (i.e., when $\rh^*(v)>0$), and that $G(v)$ is decreasing for $v\geq v_c$, \eqref{eq::uoptimal} requires the optimal reward to increase\footnote{Recall that $\frac{\partial}{\partial v}\rh^*(v)=-u^*(v)$.} at the fixed rate of $\frac{1}{\tau}$ for $v\geq v_0$, hence eventually violating the constraint $\rh^*(v)\leq p$. If $v_0<v_c$, a similar argument shows that the resulting optimal solution violates $\rh^*(v)\geq 0$ at some $v<v_0$.

This proves that $G(v)+\mu_1(v)-\mu_2(v)\leq0$ for all $v\in\mathbb{R}$.
Therefore, while off the boundary the optimal reward  decreases at the fixed rate of $\frac{1}{\tau}$. The optimal reward stays at $\hat r^*(v)=p$ for $v\leq v_L$ for some $v_L\leq v_c$ and at $\hat r^*(v)=0$ for $v\geq v_H$ for some $v_H\geq v_c$.
Note that to satisfy $G(v)+\mu_1(v)-\mu_2(v)=0$ for $v\leq v_L$ ($v\geq v_H$), $\mu_2$ ($\mu_1$) has to increase (increase)\footnote{Note the sign of $\mu_1$ vs. $\mu_2$.} at the same rate at which $G(v)$ increases (decreases) for $v\leq v_L$ ($v\geq v_H$). This also enforces $r^*_0=p$. In summary,
\begin{align}
\label{eq::uoptimal2}
u^*(v)=\begin{cases}
\frac{1}{\tau}, &v_L\leq v\leq v_H\\
0, &\text{otherwise},
\end{cases}
\end{align}
where,
\begin{align}
\label{eq::optsolution1}
G(v_L)+\mu_1^0=0,\quad G(v_H)+\mu_1^0=0,\quad v_H-v_L=p\tau.
\end{align}
Evaluating $G(v)+\mu_1(v)-\mu_2(v)=0$ at $v=-\infty$, we get 
$\mu_1^0-\mu_2^0-\lambda+\Phi(\theta-c)=0$, 
% \begin{align}
% \mu_1^0-\mu_2^0-\lambda+\Phi(\theta-c)=0,
% \end{align}
hence satisfying optimality condition for $r_0^*$ given by \eqref{eq::FOC1a}, since this zeros the coefficient of $r_0^*$ in \eqref{eq::L2}. As for the optimal cutoff, given that
\begin{align}
\frac{\partial}{\partial c}L(u^*,r_0^*,\lambda,\mu_1,\mu_2)=-c\phi(\theta-c)+\lambda(1+\frac{\partial}{\partial c}\E_{v|v_i=c}[r^*(v)]),
\end{align}
we can reach at
\begin{align}
\label{eq::optimallambda}
\lambda=\frac{c^*\phi(\theta-c^*)}{1+\normcdf{v_L-\mu_c^*}{\sigma_v}-\normcdf{v_H-\mu_c^*}{\sigma_v}}.
\end{align}
Putting \eqref{eq::optsolution1}, \eqref{eq::optimallambda}, and the indifference equation for cutoff together we can solve for the optimal reward program as given in Theorem~\ref{theorem::optimallinear}. $\hfill\blacksquare$

%%%%%%%%%%%%%%%%%%%%%%%%%%%%%%%%%%%%%%%%%%%%%%%%%%%%%%%%%%%%%%%%%%%%%%%%%%%%%
%%%%%%%%%%%%%%%%%%%%%%%%%%%%%%%%%%%%%%%%%%%%%%%%%%%%%%%%%%%%%%%%%%%%%%%%%%%%%
\noindent\textit{Proof of Lemma~\ref{lemma::Profupperbound2}.}
Choosing
$\lambda(x)=\lambda_0\times{\bf 1}\{x\leq0\}$ in \eqref{eq::Lagrangian}, we get the Lagrangian for the relaxed problem in which we drop all the constraints except for the indifference equation:
\begin{align}
\label{eq::proof::Lagrangian}
L(\rh,c,\lambda)=\E_v[(p-\rh(v)) \Phi(\frac{v-c}{\sigma_\ep})]+\lambda_0(c+\E_{v|v_i=c}[\hat r(v)]-p).
\end{align}
Recalling the identity
\begin{align}
\E_v[(p-\rh(v)) \Phi(\frac{v-c}{\sigma_\ep})]=\int_{c}^{\infty}(p-\E_{v|v_i}[\hat r(v)])\normpdf{\theta-v_i}{\sigma}dv_i,
\end{align}
we can write the
first order condition for optimal cutoff $c$ as
\begin{align}
0=-(p-\E_{v|v_i=c}[\hat r(v)])\normpdf{\theta-c}{\sigma}-\lambda_0 (-1+\frac{\tau}{\sigma_v}\E_{v|v_i=c}[(\frac{\mu_c-v}{\sigma_v})\hat r(v)]),
\end{align}
which yields
\begin{align}
\label{eq::proof::FOC1B}
\lambda_0=\frac{c\normpdf{\theta-c}{\sigma}}{1-\frac{\tau}{\sigma_v}\E_{v|v_i=c}[(\frac{\mu_c-v}{\sigma_v})\hat r(v)]}.
\end{align}
We can characterize the optimal $\rh(\cdot)$, similar to Theorem~\ref{theorem::optimal-BSpread}, by looking at the weight of $\rh(v)$ in \eqref{eq::proof::Lagrangian}, which is
$-\normcdf{v-c}{\sigma_\ep}\normpdf{v-\theta}{\sigma_\theta}+\lambda_0 \normpdf{v-\mu_c}{\sigma_v}$.
Using a similar approach, this leads to $\rh^*(v)=p\times {\bf 1}\{v<v_c\}$, where $v_c$ is the unique solution of
\be
\label{eq::proof::vc}
\frac{\lambda_0}{\normpdf{\theta-c}{\sigma}}=\frac{\normcdf{v_c-c}{\sigma_\ep}}{\normpdf{v_c-c}{\sigma_\ep}}.
\ee
The optimal expected profit for the relaxed problem is thus
\begin{align}
\E_v[(p-\rh^*(v)) \Phi(\frac{v-c}{\sigma_\ep})]=p\int_{v_c}^{\infty}\normcdf{v-c}{\sigma_\ep}\normpdf{\theta-v}{\sigma_\theta}dv.
\end{align}
The proof is now complete on noting that the optimal profit for the relaxed problem serves as an upper bound on the optimal profit for the original problem.
$\hfill\blacksquare$

%%%%%%%%%%%%%%%%%%%%%%%%%%%%%%%%%%%%%%%%%%%%%%%%%%%%%%%%%%%%%%%%%%%%%%%%%%%%%
\noindent\textit{Proof of Theorem~\ref{theorem::performance}.}
%WLOG, assume that $\sigma^2=\sigma_\ep^2+\sigma_\theta^2=1$.\\
\noindent i) When $\frac{\sigma_\ep}{\sigma_\theta}\geq\frac{p}{\sqrt{2\pi}\sigma}$, the optimal solution to the relaxed problem characterized in Lemma~\ref{lemma::Profupperbound2} coincides with the optimal spread-constrained reward program characterized in Theorem~\ref{theorem::optimallinear}, making them both optimal.

\noindent ii) We prove this by showing that the expected profit resulted from $\rh^*_{RC}(\cdot)$ approaches the upper bound $\Pi_1^H$ as $\sigma_\ep\to0$, proving the asymptotic convergence of both to the optimal expected profit.
Using the characterization of $\rh^*_{RC}(\cdot)$ in Theorem~\ref{theorem::optimallinear}, we can show that as $\sigma_\ep\to 0$,
then $c\to0$, $v_L\to 0$, $v_H\to p$. On the other hand, $\E_{v|v_i}[\rh^*_{RC}(v)]\to \rh^*_{RC}(v_i)$ as $\sigma_\ep\to 0$. As a result,
\begin{align}
\E_v[\Pi(p,\rh^*_{RC}(v))]\to\int_{0}^{\infty}(p-\rh^*_{RC}(v_i))\normpdf{\theta-v_i}{\sigma}dv_i\to\int_{0}^{p}v_i\normpdf{\theta-v_i}{\sigma}dv_i+p\normcdf{\theta-p}{\sigma},\nonumber\\
\end{align}
which is the same as the upper bound $\Pi_1^H$ given by \eqref{eq::profubound1}.
$\hfill\blacksquare$